\documentclass[a4paper,11pt]{article}
\usepackage[top=2.5cm,bottom=2.5cm,left=2.2cm,right=2.2cm]{geometry}
\usepackage{amsfonts}
\usepackage{mathrsfs,amscd,amssymb,amsthm,amsmath,bm,graphicx,psfrag,subfigure,url,mathtools}
\usepackage{pict2e}
\usepackage{stfloats}
\usepackage{psfrag,amsmath}
\usepackage{tikz}
\usepackage{indentfirst}
\usepackage{hyperref}
\usepackage{bookmark}
\usepackage{enumerate}
\usepackage{latexsym,euscript,epic,eepic,color}
\usepackage{multirow}
\usepackage{multicol}
\usepackage{longtable}
\usepackage{adjustbox}
\usepackage[all]{xy}
\usepackage{setspace}
\allowdisplaybreaks
\usepackage{authblk}
\usepackage{amsfonts}
\usepackage{amsmath}
\usepackage{mathrsfs,color}
\usepackage{epstopdf}
\usepackage{graphicx}
\usepackage{tikz}
\usetikzlibrary{positioning}
\usetikzlibrary{fadings}
\usetikzlibrary{patterns}
\usetikzlibrary{calc}
\usetikzlibrary{shadings}
\pgfsetlayers{background,main,foreground}

\makeatletter

\renewcommand{\@seccntformat}[1]{{\csname the#1\endcsname}{\normalsize .}\hspace{.5em}}
\makeatother

\usepackage{ifpdf}
\usepackage{indentfirst}

\def \[{\begin{equation}}
\def \]{\end{equation}}

\newtheorem{thm}{Theorem}[section]

\newtheorem{claim}{Claim}
\newtheorem{lem}[thm]{Lemma}
\newtheorem{cor}[thm]{Corollary}

\newtheorem{pb}{Problem}

\newtheorem{conj}[thm]{Conjecture}

\begin{document}
\baselineskip=0.23in
\title{\bf Spectral extrema of graphs with fixed size: forbidden a fan graph,  friendship graph  or theta graph\thanks{{\it Email addresses}: lscmath@ccnu.edu.cn (S.C. Li),\ 15537204105@163.com (S.S. Zhao),\ zoulantao123@163.com (L.T. Zou)}}
\author[,1,2]{Shuchao Li}
\author[1]{Sishu Zhao}
\author[,3]{Lantao Zou}
\affil[1]{School of Mathematics and Statistics, and Hubei Key Lab--Math. Sci.,\linebreak Central China Normal University, Wuhan 430079, China}
\affil[2]{Key Laboratory of Nonlinear Analysis \& Applications (Ministry of Education),\linebreak Central China Normal University, Wuhan 430079, China}
\affil[3]{School of Mathematics and Statistics, Central South University, Changsha 410083, China}
\date{\today}
\maketitle
\begin{abstract}
It is well-known that the Brualdi-Hoffman-Tur\'an-type problem inquiries about the maximum spectral radius \( \lambda(G) \) of an \( F \)-free graph \( G \) with \( m \) edges. This can be regarded as a spectral characterization of the existence of the subgraph \( F \) within \( G \). A significant contribution to this problem was made by Nikiforov (2002). He proved that for every \( K_{r + 1} \)-free graph with \( m \) edges, \( \lambda(G) \leqslant \sqrt{2m(1 - 1/r)} \). Let \( \theta_{1,p,q} \) denote the theta graph, which is constructed by connecting two vertices with 3 internally disjoint paths of lengths 1, \( p \), and \( q \) respectively. Let \( F_k \) be the fan graph, that is, the join of a \( K_1 \) and a path \( P_{k - 1} \). Let \( F_{k,3} \) be the friendship graph, obtained by having \( k \) triangles share a common vertex. In this paper, we utilize the \( k \)-core method and spectral techniques to address some spectral extrema of graphs with a fixed number of edges. Firstly, we demonstrate that for \( m \geqslant \frac{9}{4}k^6 + 6k^5 + 46k^4 + 56k^3 + 196k^2 \) and \( k \geqslant 3 \), if \( G \) is \( F_{2k + 2} \)-free, then \( \lambda(G) \leqslant \frac{k - 1 + \sqrt{4m - k^2 + 1}}{2} \). Equality holds if and only if \( G \cong K_k \vee (\frac{m}{k}-\frac{k - 1}{2})K_1 \). This validates a conjecture by Yu, Li, and Peng [Discrete Math. 348 (2025) 114391] and refines a recent result by Li, Zhai, and Shu [European J. Combin. 120 (2024) 103966]. Secondly, we show that for \( m \geqslant \frac{9}{4}k^6 + 6k^5 + 46k^4 + 56k^3 + 196k^2 \) with \( k \geqslant 3 \), if \( G \) is \( F_{k,3} \)-free and has \( m \) edges, then \( \lambda(G) \leqslant \frac{k - 1 + \sqrt{4m - k^2 + 1}}{2} \). Equality holds precisely when \( G \cong K_k \vee (\frac{m}{k}-\frac{k - 1}{2})K_1 \). This confirms a conjecture put forward by Li, Lu, and Peng [Discrete Math. 346(2023)113680]. Finally, we identify the \( \theta_{1,p,q} \)-free graph with \( m \) edges that possesses the largest spectral radius, where \( q \geqslant p \geqslant 3 \) and \( p + q \geqslant 2k + 1 \). A further research problem is also proposed. 

\vskip 0.2cm
\noindent {\bf Keywords:}  Spectral radius; Friendship graph; Fan graph; Theta graph; Extremal graph\vspace{1mm}

\noindent {\bf AMS Subject Classification:} 05C50;\, 05C35
\end{abstract}

\section{\normalsize Introduction}
In this paper, we focus on simple and finite graphs. Unless specified otherwise, we adhere to traditional notation and terminology (see, for example, Godsil and Royle~\cite{GR-01}, West~\cite{WST-96}).  

Let \( G = (V(G), E(G)) \) be a graph with vertex set \( V(G) = \{v_1, \dots, v_n\} \) and edge set \( E(G) = \{e_1, \dots, e_m\} \), where \( n \) (\textit{order}) and \( m \) (\textit{size}) denote the number of vertices and edges of \( G \), respectively. Let \( A(G) \) be the adjacency matrix of \( G \). Since \( A(G) \) is real symmetric, its eigenvalues are real and can be ordered as \( \lambda_1(G) \geqslant \dots \geqslant \lambda_n(G) \). The \textit{spectral radius} \( \lambda(G) \) of \( G \) is defined as \( \max\{|\lambda_1(G)|, \dots, |\lambda_n(G)|\} \). By the Perron-Frobenius theorem, \( \lambda(G) = \lambda_1(G) \).  
For two graphs \( G \) and \( H \), define \( G \cup H \) as their disjoint union (no shared vertices). Their join \( G \vee H \) is constructed from \( G \cup H \) by adding all possible edges between vertices of \( G \) and vertices of \( H \).  

In 1985, Brualdi and Hoffman~\cite{BH-85} initiated the problem on characterizing graphs of given size having maximal spectral radius.  In particular, they posed the following conjecture.
\begin{conj}\label{conj-1}
Let $G$ be a graph of size $m$ without isolated vertices.  If $m ={a \choose 2}+b$  with $0 \leqslant b <a,$  then $\lambda(G)\leqslant \lambda(K_b\vee (K_{a-b}\cup K_1)),$  with equality if and only if $G\cong K_b\vee (K_{a-b}\cup K_1).$
\end{conj}
Some special cases of Conjecture~\ref{conj-1} were confirmed by Brualdi and Hoffman~\cite{BH-85}, Friedland~\cite{F-85} and
Stanley~\cite{S-87}.  Conjecture~\ref{conj-1} was fully resolved by Rowlinson~\cite{R-88}.

As an analogue of Conjecture~\ref{conj-1}, Bhattacharya, Friedland and Peled~\cite{BFP-08} posed the following conjecture.
\begin{conj}\label{conj-2}
Let $G$ be a bipartite graph of size $m$ with bipartite sets $S$  and $T,$  where $2\leqslant |S|\leqslant |T|$ and $0<m<|S||T|.$  If $G$ achieves the maximum spectral radius, then $G$ is obtained from a complete bipartite graph by adding one vertex and a corresponding number of edges.
\end{conj}
Conjecture~\ref{conj-2} has been confirmed for certain special cases by Bhattacharya, Friedland and Peled~\cite{BFP-08}, Chen et al.~\cite{CFK-10}, Das et al.~\cite{DC-13} and Liu and Weng [20]. For further developments regarding Conjecture~\ref{conj-2},  we direct the reader to \cite{YGL-24,ZLZ-22}.  

Let $H$ be a graph. A graph $G$ is said to be $H$-\textit{free}, if it does not contain $H$ as a subgraph. Let $\mathcal{G}(m,H)$ be the set of all $H$-free graphs with  size $m$. The following is the well-known Brualdi-Hoffman-Tur\'an-type problem, which has been drawing increasing attention from researchers \cite{LY2024,C01}.
\begin{pb}[Brualdi-Hoffman-Tur\'an-type problem]\label{pb-1}
  What is the maximum spectral radius $\lambda(G)$ of an $H$-free graph $G$ with $m$ edges$?$
\end{pb}
Generally speaking, the study of Problem ~\ref{pb-1} has a close connection with triangles. Note that a triangle can be regarded as either a cycle \( C_3 \) or a complete graph \( K_3 \). The research on Problem ~\ref{pb-1} originates from the study of triangles. In 1970, Nosal~\cite{N-70} proved that \( \lambda(G) \leqslant \sqrt{m} \) for every graph \( G \) in \( \mathcal{G}(m, K_3) \). Lin, Ning, and Wu \cite{LNW-21} extended Nosal’s result, showing that for a non-bipartite \( C_3 \)-free graph \( G \) with \( m \) edges, \( \lambda(G) \leqslant \sqrt{m - 1} \), and equality holds if and only if \( G \cong C_5 \). Under the same condition, Zhai and Shu \cite{ZS-22} improved Lin, Ning, and Wu’s results, demonstrating that \( \lambda(G) \leq \lambda(SK_{2, \frac{m - 1}{2}}) \), with equality if and only if \( G \cong SK_{2, \frac{m - 1}{2}} \). Here, \( SK_{2, \frac{m - 1}{2}} \) is obtained by subdividing an edge of \( K_{2, \frac{m - 1}{2}} \).

The second aspect of Problem \ref{pb-1} involves the family of odd cycles that include \( C_3 \). Sun and Li \cite{SL-23} showed that if \( G \) is a non-bipartite \( \{C_3, C_5\} \)-free graph with \( m \) edges, then \( \lambda(G) \leqslant \sqrt[4]{\sum_{u \in V_G} d_u^2 - m + 4q + 5} \), and equality holds if and only if \( G \cong C_7 \), where \( q \) denotes the number of 4-cycles in \( G \). Let \( r(m) \) be the largest root of the equation \( x^4 - x^3 - (m - 3)x^2 + (m - 4)x + m - 5 = 0 \). Li, Peng \cite{LP-2022}, and Sun, Li \cite{SL-23} independently investigated the further stability result as follows: Let \( G \) be a non-bipartite \( \{C_3, C_5\} \)-free graph with \( m \) edges; then \( \lambda(G) \leqslant r(m) \), and equality holds if and only if \( G \cong RK_{2, \frac{m - 3}{2}} \) when \( m \) is odd. Here, \( RK_{2, \frac{m - 3}{2}} \) is obtained by replacing one edge of the complete bipartite graph \( K_{2, \frac{m - 3}{2}} \) with \( P_5 \). The case when \( m \) is even was recently solved by Li and Yu \cite{LY2024}.  

Li, Sun, and Yu \cite{LSY01} also demonstrated that, for a \( \{C_3, C_5, \ldots, C_{2k + 1}\} \)-free graph \( G \), the inequality \( \lambda_1^{2k} + \lambda_2^{2k} \leq \frac{\text{Tr}(A^{2k}(G))}{2} \) holds, where \( \text{Tr}(\cdot) \) denotes the trace of the corresponding matrix. All the relevant extremal graphs have been characterized. Furthermore, if \( G \) is non - bipartite, then  

\[\notag
\lambda^{2k}(G) \leqslant \frac{\text{Tr}(A^{2k}(G))}{2} - \left(2\cos\frac{\pi}{k + 2}\right)^{2k}
\]  
Equality is achieved if and only if \( k = 1 \) and \( G \cong C_5 \). Evidently, when \( k = 1 \), this result coincides with the earlier finding by Lin, Ning, and Wu \cite{LNW-21} mentioned above.  

The third aspect of Problem \ref{pb-1} involves two scenarios: either identifying an edge of a triangle with an edge of another cycle, or having a vertex of a triangle shared with vertices of other cycles. Let \( \theta_{t,p,q} \) denote the theta graph, constructed by connecting two vertices via 3 internally disjoint paths of lengths \( t \), \( p \), and \( q \). Sun, Li, and Wei \cite{SLW-23} established tight upper bounds for \( \lambda(G) \) when \( G \) belongs to \( \mathcal{G}(m, \theta_{1,2,3}) \) and \( \mathcal{G}(m, \theta_{1,2,4}) \), respectively. From these, one can deduce the graph within \( \mathcal{G}(m, C_5) \) or \( \mathcal{G}(m, C_6) \) that attains the largest spectral radius (refer to \cite[Theorem 1.2]{C95}). Recently, Lu, Lu, and Li \cite{LLL24} identified the graph in \( \mathcal{G}(m, \theta_{1,2,5}) \) with the maximum spectral radius. 

Let \( B_{r + 1} \) represent the \( (r + 1) \)-book graph, formed by \( r + 1 \) triangles sharing a common edge. Nikiforov \cite{N-21} determined the graph in \( \mathcal{G}(m, B_{r + 1}) \) that has the largest spectral radius.  

Recently, Li, Zhai and Shu \cite{C01} obtained Theorem~\ref{cj1-4}, which confirms a conjecture proposed by Li for $m=\Omega(k^4)$ (see also~\cite{LLH2022}). 
\begin{thm}[\cite{C01}]\label{cj1-4}
Let $k\geqslant 3$ and $m=\Omega(k^4)$. 
If $G\in  \mathcal{G}(m, \theta_{1,2,2k-1})\cup\mathcal{G}(m, \theta_{1,2,2k})$, then $\lambda(G)\leqslant \frac{k-1+\sqrt{4m-k^2+1}}{2}$  with equality if and only if $G\cong K_k\vee (\frac{m}{k}-\frac{k-1}{2})K_1.$
\end{thm}
In fact, Theorem \ref{cj1-4} can be used to deduce the following theorem, which also confirms a conjecture put forward by Zhai, Lin, and Shu \cite{C95}. 

\begin{thm}[\cite{C01}]\label{cj1-3}
Let $k$ be a fixed positive integer and $G$ be a graph of sufficiently  large size $m$ without isolated vertices. If $\lambda(G)\geqslant \frac{k-1+\sqrt{4m-k^2+1}}{2},$ then $G$ contains a cycle $C_t$ for every $t\leqslant 2k+2$, unless $G\cong K_k\vee (\frac{m}{k}-\frac{k-1}{2})K_1.$
\end{thm}
Motivated by Theorem~\ref{cj1-4}, it is natural to consider the following problem.
\begin{pb}\label{pb-2}
 What is the maximum spectral radius of graphs among $\mathcal{G}(m,\theta_{1,p,q})$ for $q \geqslant p\geqslant 3?$
\end{pb}
Recall that \( F_{k,3} \) is the friendship graph formed by \( k \) triangles sharing a common vertex. Li, Lu, and Peng \cite{LLP-23} demonstrated that, for graphs in \( \mathcal{G}(m, F_{2,3}) \), the maximum spectral radius is \( \frac{1 + \sqrt{4m - 3}}{2} \), and the corresponding extremal graph is \( K_2 \vee \frac{m - 1}{2}K_1 \). Additionally, they proposed the following conjecture. 
\begin{conj}[\cite{LLP-23}]\label{cj1-5}
Let \( k \geqslant 3 \) be a fixed integer and \( m \) be sufficiently large. If \( G \in \mathcal{G}(m, F_{k,3}),\) then  
$
\lambda(G) \leqslant \frac{k - 1 + \sqrt{4m - k^2 + 1}}{2}
$  
with equality if and only if \( G \cong K_k \vee \left( \frac{m}{k} - \frac{k - 1}{2} \right) K_1 \).
\end{conj}
Let \( F_k = K_1 \vee P_{k - 1} \) denote the fan graph on \( k \) vertices, where the vertex with degree \( k - 1 \) is referred to as the \textit{central vertex}. Recently, Yu, Li, and Peng \cite{YLP2024}  put forward the following conjecture.
\begin{conj}[\cite{YLP2024}]\label{conj-6}
Let $k \geqslant 2$ be ﬁxed and $m$ be suﬃciently large. If $G\in  \mathcal{G}(m, F_{2k+1})$ or $ G\in\mathcal{G}(m, F_{2k+2}),$ then
$
\lambda(G)\leqslant\frac{k-1+\sqrt{4m-k^2+1}}{2}
$
with equality if and only if $G \cong K_k\vee(\frac{m}{k}-\frac{k-1}{2})K_1$.
\end{conj}
Yu, Li, Peng \cite{YLP2024}, and Zhang, Wang \cite{zw2024} independently investigated the case of \( k = 2 \) for \( F_{2k + 1} \)-free graphs, while Gao and Li \cite{GL-2024} looked into the case of \( k = 2 \) for \( F_{2k + 2} \)-free graphs.  

Motivated by \cite{C01,YLP2024,zw2024}, in this paper, we consider \( F_{2k + 2} \)-free graphs and \( F_{k, 3} \)-free graphs respectively for \( k \geqslant 3 \). We will present a unified approach to address Problem~\ref{pb-2}, Conjectures~\ref{cj1-5} and~\ref{conj-6}.

Our first main result determines the largest spectral radius of an \( F_{2k + 2} \)-free graph with size \( m \) for \( k \geqslant 3 \), and characterizes the corresponding extremal graph.  
\begin{thm}\label{thm2.1}
Let $k\geqslant3$ and $m\geqslant \frac{9}{4}k^6+6k^5+46k^4+56k^3+196k^2$. If $G\in\mathcal{G}(m, F_{2k+2}),$ then
$
\lambda(G)\leqslant\frac{k-1+\sqrt{4m-k^2+1}}{2}
$
with equality if and only if $G \cong K_k\vee(\frac{m}{k}-\frac{k-1}{2})K_1$.
\end{thm}
Observe that every $F_{2k+1}$-free must be $F_{2k+2}$-free, that is, $\mathcal{G}(m,F_{2k+1})\subseteq \mathcal{G}(m,F_{2k+2})$, and $K_{k}\vee (\frac{m}{k}-\frac{k-1}{2})K_1\in\mathcal{G}(m,F_{2k+1}) $. Hence, Conjecture~\ref{conj-6} is a direct consequence of Theorem~\ref{thm2.1}.

Our second main result determines the largest spectral radius of $F_{k,3}$-free graph of size $m$ for $k\geqslant 3$, and identifies the corresponding extremal graph. By Corollary~\ref{thm1.3} below, Conjecture~\ref{cj1-5} follows immediately. 
\begin{cor}\label{thm1.3}
Let $k\geqslant3$ and 
 $m\geqslant \frac{9}{4}k^6+6k^5+46k^4+56k^3+196k^2$.
If $G\in\mathcal{G}(m,F_{k,3}),$ then
$
\lambda(G)\leqslant\frac{k-1+\sqrt{4m-k^2+1}}{2}
$
with equality if and only if $G\cong K_k\vee(\frac{m}{k}-\frac{k-1}{2})K_1$.
\end{cor}

Our last main result determines the largest spectral radius of $\theta_{1,p,q}$-free graph of size $m$ for $q\geqslant p\geqslant 3$, and also characterizes the corresponding extremal graph. Consequently, Corollary~\ref{thm2.2} below resolves Problem~\ref{pb-2} for $p+q\geqslant 7$.
\begin{cor}\label{thm2.2}
Let $k\geqslant 3$ and $m\geqslant \frac{9}{4}k^6+6k^5+46k^4+56k^3+196k^2$. If $G\in\mathcal{G}(m, \theta_{1,p,q})$ or $G\in \mathcal{G}(m, \theta_{1,r,s})$ with $q\geqslant p\geqslant 3,s\geqslant r\geqslant 3, p+q=2k+1$ and $r+s=2k+2,$ then
$
\lambda(G)\leqslant\frac{k-1+\sqrt{4m-k^2+1}}{2}
$
with equality if and only if $G \cong K_k\vee(\frac{m}{k}-\frac{k-1}{2})K_1$.
\end{cor}

\noindent{\bf Organization.} We start by introducing relevant notation and presenting preliminary results in Section~\ref{sec-2}.  
In Section~\ref{sec-3}, we conduct an in-depth characterization of the local structure of the extremal graph, which serves as the theoretical foundation for proving our main results. In Section~\ref{sec-4}, we provide proofs for Theorem~\ref{thm2.1}, as well as Corollaries~\ref{thm1.3} and \ref{thm2.2}. Some concluding remarks are offered in the final section.
\section{\normalsize Preliminaries} \setcounter{equation}{0}\label{sec-2}
\noindent{\bf Notation.}  
For a graph \( G \) and a vertex \( u \in V(G) \), let \( N_G(u) \) denote the neighborhood of \( u \) in \( G \), and \( N_G[u] = N_G(u) \cup \{u\} \). Let \( e(G) = |E(G)| \) (the \textit{size} of \( G \)) and \( |G| = |V(G)| \) (the \textit{order} of \( G \)).  For two vertex-disjoint subsets \( S, T \subseteq V(G) \), let \( E_G(S, T) \) be the set of edges with one endpoint in \( S \) and the other in \( T \), and \( e_G(S, T) = |E_G(S, T)| \). Let \( G[S] \) denote the subgraph induced by \( S \), and \( E_G(S) \) the edge set of \( G[S] \); we also write \( e_G(S) = |E_G(S)| \).  For vertex subsets \( S, T \subseteq V(G) \), for simplicity, we use \( N_S(T) \) to denote \( \bigcup_{u \in T} N(u) \cap S \) and \( N_S[T] = N_S(T) \cup T \). Additionally, we identify \( S \) (resp. \( T \)) with the induced subgraph \( G[S] \) (resp. \( G[T] \)) when context permits.  Subscripts (e.g., \( N_G(u) \), \( E_G(S, T) \)) may be omitted if their meaning is clear from context. 

Throughout this text, \( P_n \), \( C_n \), and \( K_n \) denote the path, cycle, and complete graph on \( n \) vertices, respectively. Let \( G - uv \) represent the graph derived from \( G \) by removing the edge \( uv \in E(G) \). Let \( G + uv \) represent the graph formed from \( G \) by adding the edge \( uv \notin E(G) \); this notation extends naturally to cases involving the addition or removal of more than one edge. The symbol \( \sim \) indicates that the two vertices in question are adjacent.

Since the adjacency matrix $A(G)$ is irreducible and nonnegative for a connected graph, by Perron-Frobenius theorem, we know that the largest eigenvalue of $A(G)$ is equal to the spectral radius $\lambda(G)$ of $G$, and there exists a positive eigenvector $x$ of $A(G)$ corresponding to $\lambda(G)$. The eigenvector $x$ mentioned above is the \textit{Perron vector} of $G$ with coordinate $x_v$ corresponding to the vertex $v\in V(G)$. A vertex $u^*$ is said to be an \textit{extremal vertex} if $x_{u^*}=\max_{u\in V(G)}\{x_{u}\}$.

Let \( S_{n,k} = K_k \vee (n - k)K_1 \). Let \( S_{n,k}^+ \) denote the graph formed by adding an edge within the independent set of \( S_{n,k} \). Additionally, let \( M_t \) be the graph with \( t \) vertices, consisting of a matching containing \( \left\lfloor \frac{t}{2} \right\rfloor \) edges, along with one additional vertex if \( t \) is odd.  \vspace{3mm}

\noindent{\bf Some basic lemmas.} 
\begin{lem}[\cite{ZS-23}]\label{lem5.1}
Let $A$ and $ A'$ be the adjacency matrices of two connected graphs $G$ and $G'$ with the same vertex set. Suppose that $N_G(u) \subsetneqq N_{G'}(u)$ for some vertex $u$. If the Perron vector $\boldsymbol{x}$ of $G$ satisfies $\boldsymbol{x}^TA^{\prime }\boldsymbol{x}\geqslant \boldsymbol{x}^TA\boldsymbol{x}$, then $ \lambda(G^{\prime})>\lambda(G).$
\end{lem}

\begin{lem}[\cite{C95}]\label{lem2.4a}
Let $F$ be a $2$-connected graph and $G$ attain the maximum spectral radius in $\mathcal G(m,F)$. Then $G$ is connected. Moreover, if $u^*$ is an extremal vertex of $G$, then there exists no cut vertex in $G[V(G)\setminus\{u^*\}]$.
\end{lem}
Note that given a graph $G$ and  a vertex $u\in V(G)$, if  $G[N(u)]$ contains a $P_{2k+1}$, then we can find an $F_{2k+2}$ in $G$. Hence, the following result follows immediately.
\begin{lem}\label{lem2.5a}
Let $G$ be a graph in  $\mathcal{G}(m,F_{2k+2})$. Then for all $u\in V(G)$, the graph $G[N(u)]$ is $P_{2k+1}$-free.  
\end{lem}

In the following, let $G^*$ be the  graph in  $\mathcal{G}(m,F_{2k+2})$ having the maximum spectral radius. By Lemma~\ref{lem2.4a}, one can see that $G^*$ is connected.  By Perron-Frobenius theorem, there exists a positive eigenvector $\boldsymbol{x}$ corresponding to $\lambda(G^*)$ with coordinate $x_v$ corresponding to the vertex $v\in V(G)$. 
We may assume that $\lambda(G^*)=\lambda$ and $x_{u^*}=\max_{u\in V(G^*)}x_u=1$ for some $u^*\in V(G^*)$. Furthermore, we denote $R=N(u^*)$, $S=V(G^*)\backslash N[u^*]$ and $d_R(u)=|N_R(u)|$ for each vertex $u\in V(G^*)$.

Let $\gamma=-\frac{k(k-1)}{2}$. 
Note that $K_k\vee(\frac{m}{k}-\frac{k-1}{2})K_1\in \mathcal{G}(m,F_{2k+2}).$ It follows that $\lambda\geqslant \lambda(K_k\vee(\frac{m}{k}-\frac{k-1}{2})K_1)$. So from then on we may assume that
 \begin{eqnarray}\label{1b}
     \lambda^2-(k-1)\lambda\geqslant m+\gamma.
   \end{eqnarray}
  Note that $\lambda=\lambda x_{u^*}=\sum_{u\in R}x_u$, and
\begin{align*}
 \lambda^2=\lambda^2x_{u^*}=d_{G^*}(u^*)x_{u^*}+\sum_{u\in R}d_R(u)x_u+\sum_{w\in S}d_R(w)x_w.
\end{align*}
Hence,  
\begin{eqnarray}\label{3b}
 \lambda^2-(k-1)\lambda=d_{G^*}(u^{*})+\sum_{u\in R}(d_R(u)-k+1)x_u+\sum_{w\in S}d_R(w)x_w.
\end{eqnarray}

For an arbitrary subset $L$ of $R$, we define
\begin{eqnarray}\label{4}
\eta(L)=\sum_{u\in L}(d_L(u)-k+1)x_u-e(L).
\end{eqnarray}
In particularly, if $L=\emptyset$, we define $\eta(L)=0$. 
Together with \eqref{3b} and \eqref{4}, we get
    \begin{align}
   \lambda^2-(k-1)\lambda&=d_{G^*}(u^*)+\eta(R)+e(R)+\sum_{w\in S}d_R(w)x_w\label{5.0b}\\
   &\leqslant  d_{G^*}(u^*)+\eta(R)+e(R)+e(R,S)\label{5b}\\
   &=\eta(R)+m-e(S).\notag 
   \end{align}
Combining \eqref{1b} and \eqref{5b} gives us
 \begin{eqnarray}\label{6b}
   \eta(R)\geqslant e(S)+\gamma\geqslant\gamma.
   \end{eqnarray}
By \eqref{5b} and \eqref{6b}, one sees if $\eta(R)=\gamma$, then $e(S)=0$ and $x_w=1$ for each $w\in S$.

We now recall the terminology of $k$-core, introduced by Seidman~\cite{SBS-1983} in 1983. A \textit{$k$-core} of a graph $G$ is the largest induced subgraph of $G$ such that its minimum degree is at least $k$. It is obvious that a $k$-core can be obtained iteratively from $G$ by deleting the vertices of degree at most $k-1$ until the resulting graph is empty or is of minimum degree at least $k$. It is known that $k$-core is well-defined, that is, it does not depend on the order of vertex deletion. A graph is referred to as $(k-1)$-\textit{degenetate} if its $k$-core is empty. It brings a breakthrough in extremal graph theory (see \cite{AKB-2003} and \cite{L-2017} for details). Nikiforov~\cite{VN-2014} was the first to utilize these notions to study spectral extremal graph theory. Now the core of graph is a key tool, which was used to study the spectral graph theory (see \cite{HFD-2023,C01}).

In the following, we shall introduce a variable on $R$ as follows. Now we denote by $L^c$ the vertex set of the $(k-1)$-core of $G^*[L]$. It is obvious that $L^c\subseteq L$ for every subset $L$ of $R$. And if $L=\emptyset$, then $L^c=\emptyset$. We need the following lemmas.
\begin{lem}[\cite{C01}]\label{lem3.2a}
For every subset $L$ of $R$, we have $\eta(L)\leqslant \eta(L^c)$ with equality if and only if $L=L^c$.
\end{lem}


Let $\mathcal J$ be the family of connected components in $G^*[R^c]$ and $|\mathcal J|$ be the number of members in $\mathcal J$. By the definition of $(k-1)$-core, we have $\delta(J)\geqslant k-1$ for each $J\in\mathcal J$. Therefore, it follows from \eqref{4} that for each $J\in\mathcal J,$
\begin{eqnarray}\label{7}
   \eta(V(J))&\leqslant& \sum_{u\in V(J)}(d_J(u)-k+1)-e(J)\\
   &=&e(J)-(k-1)|J|. \notag
\end{eqnarray}
Equality in \eqref{7} holds if and only if $x_u=1$ for each $u\in V(J)$ with $d_J(u)\geqslant k$. 
We denote by $\mathcal{L}_{|J|,k-1}$ the family of graphs obtained from $S_{|J|,k-1}^+$ by deleting an arbitrary edge. 
By Lemma~\ref{lem2.5a},  $J\subseteq G^*[R]$ is $P_{2k+1}$-free, then the following lemmas in \cite{C01} still hold.

\begin{lem}[\cite{C01}]\label{lem3.3a}
Let $\mathcal J_1=\{J\in\mathcal J||J|\geqslant2k+1\}$. Then for each $J\in\mathcal J_1$,
$$
\eta(V(J))\leqslant 
\left\{
  \begin{array}{ll}
    \gamma+1, & \hbox{if $J\cong S_{|J|,k-1}^+;$} \\[5pt]
    \gamma, & \hbox{if $J\in \mathcal{L}_{|J|,k-1};$} \\[5pt]
    \gamma-\frac{1}{2}, & \hbox{otherwise.}
  \end{array}
\right.
$$
If $J\in \mathcal L_{|J|,k-1}$ and $\eta(V(J))=\gamma$, then $x_u=1$ for each $u\in V(J)$ with $d_{J}(u)\geqslant k$.
\end{lem}

In what follows, we consider the members in $\mathcal J\setminus \mathcal J_1$. Recall that $\delta(J)\geqslant k-1$ for each $J\in\mathcal{J}$. Hence, for each $J\in\mathcal{J}\setminus \mathcal J_1$, we have $k\leqslant |J|\leqslant 2k$. Now, let $\mathcal J_2$ be the subfamily of $\mathcal{J}\setminus \mathcal J_1$, in which every member does not contain any cycle of length large than $2k-2$. 
\begin{lem}[\cite{C01}]\label{lem3.4a}
For every member $J\in\mathcal J_2$, we have $\eta(V(J))\leqslant -(k-1)$.
\end{lem}
\section{\normalsize Characterizing $\mathcal{J}_1, \mathcal{J}_2, \mathcal{J}_3,\mathcal{J}_4$ and $\mathcal{J}_5$}\label{sec-3}
Recall that $\mathcal J_1=\{J\in\mathcal J||J|\geqslant2k+1\}$ and $\mathcal J_2$ is the subfamily of $\mathcal{J}\setminus \mathcal J_1$, in which every member does not contain any cycle of length larger than $2k-2$. So we let $\mathcal{J}_3$ be the subfamily of $\mathcal{J}\setminus(\mathcal{J}_1\cup\mathcal{J}_2)$, in which every member does not contain any cycle of length larger than $2k-1$. Then for every member $\hat{J}\in\mathcal{J}_3$, it contains a longest cycle of length $2k-1$ and $2k-1\leqslant|\hat{J}|\leqslant 2k$. 
Let $\mathcal{J}_4=\mathcal{J}\setminus(\mathcal{J}_1\cup\mathcal{J}_2\cup\mathcal{J}_3)$, that is, $J$ contains a longest cycle of length $2k$ 
for each $J\in \mathcal J_4$. Thus, $|J|=2k$ and $N_R(u)\subseteq V(J)$ for each $u\in V(J)$, otherwise, we obtain a $P_{2k+1}$ in $G^*[R]$, a contradiction. Moreover, let $\mathcal{J}_5$ be the subfamily of $\mathcal{J}_4$, in which $\eta(V(J))>0$ for each $J\in\mathcal J_5$. 

In this section, we firstly establish an upper bound on $\eta(V(J))$ for $J\in \mathcal J_3\cup \mathcal J_4$. Then we determine the cardinalities, respectively, for $\mathcal J_1, \mathcal J_2, \mathcal J_3$ and $\mathcal J_4$. In the remaining of our context, the notation $k$ is always referred to parameter in the forbidden graph $F_{2k+2}$.

\subsection{\normalsize Upper bound on $\eta(V(J))$ for $J\in \mathcal J_3\cup \mathcal J_4$.}

In this subsection, we establish an upper bound on  $\eta(V(J))$ for $J\in \mathcal J_3\cup \mathcal J_4$. We need the following lemma.
\begin{lem}\label{lem5.2b}
Let $s\geqslant 2$ and $G=K_1\vee H$, where graph $H$ is obtained from $K_{2s-1}$ by deleting arbitrary $s$ edges. Then for any $v\in V(H)$, there exists a $P_{2s}$ in $G$ starting from $v$.
\end{lem}
\begin{proof}
    We prove the result by induction on $s$. Obviously, the result  is true for $s=2$.  Now let $\ell\geqslant 3$ and assume the result is true for $s\leqslant \ell-1$. Let $G=K_1\vee H$, where graph $H$ is obtained from $K_{2\ell-1}$ by deleting arbitrary $\ell$ edges. For any $v\in V(H)$, we can find a $u\in N_H(v)$ such that $H-v-u$ is obtained from $K_{2\ell-3}$ by deleting at most $\ell-1$ edges. By induction, for any $w\in N_H(u)\cap (V(H)\setminus \{u,v\}),$ there exists a $P_{2\ell-2}$ in $G-\{u,v\}$  starting from $w$ and then we find a $P_{2\ell}$ in $G$ starting from $v$. This completes the proof.
\end{proof}

In order to establish an upper bound on  $\eta(V(J))$ for $J\in \mathcal J_3\cup \mathcal J_4$, we need the following key lemma.
\begin{lem} \label{lem 5*}
$\mathcal{J}_5$ is empty.
\end{lem}
\begin{proof}
Suppose to the contrary that $|\mathcal{J}_5\mathcal|\geqslant 1$. We need the following claims to complete the proof. 
\begin{claim}\label{1^}
    $e(J)>2k(k-1)$ for each $J\in \mathcal J_5$.
\end{claim}
\begin{proof}[\bf Proof of Claim~\ref{1^}]
Suppose to the contrary that $e(J)\leqslant 2k(k-1)$ for some $J\in\mathcal{J}_5$.
By  \eqref{7},
we obtain $\eta(V(J))\leqslant e(J)-(k-1)|J|\leqslant 0$, contradicting the definition of $\mathcal{J}_5$.
\end{proof}

\begin{claim}\label{2^}
$\sum_{v\in V(J)} x_v>2k-2$ for each $J\in \mathcal{J}_5$.
\end{claim}
\begin{proof}[\bf Proof of Claim~\ref{2^}] 
Suppose that there exists a $J \in \mathcal J_5$ such that $\sum_{v\in V(J)} x_v\leqslant 2k-2$. By   \eqref{4} and Claim \ref{1^}, we obtain $\eta(V(J))\leqslant(\Delta(J)-k+1)\sum_{v\in V(J)} x_v-e(J)$ $<k(2k-2)-2k(k-1)=0$, a contradiction. 
\end{proof}
\begin{claim}\label{3^}
$|\mathcal J_5| \leqslant \frac{\lambda}{2k-2}+1$.
\end{claim}
\begin{proof}[\bf Proof of Claim~\ref{3^}]
Suppose to the contrary that $|\mathcal J_5|>\frac{\lambda}{2k-2}+1$. 
Note that for each $J\in \mathcal{J}_5$ and  each $u\in V(J)$,  $N_R(u)\subseteq V(J)$, and so $d_R(u)=d_J(u)$. Next we show
\[\label{eq:3.005}
\text{$\eta(V(J))\leqslant k$\ \ for each $J\in \mathcal{J}_4.$}
\]
In fact, one may see that $|J|=2k$ for each $J\in \mathcal J_4$. Then $\Delta(J)\leqslant 2k-1$ and so $e(J)\leqslant \frac{2k-1}{2}|J|$. In view of ~\eqref{7}, we have $\eta(V(J))\leqslant e(J)-(k-1)|J|\leqslant \frac{1}{2}|J|= k$, as desired.

Recall that $e(J)\leqslant k(2k-1)$ for each $J\in \mathcal J_5$. By Claim~\ref{2^} and \eqref{eq:3.005}, for each $J\in\mathcal{J}_5$, we obtain
\begin{align*}
 (\lambda-k+1)(2k-2)
&<(\lambda-k+1) \sum_{v\in V(J)} x_v \\
& =\sum_{v\in V(J)}(x_{u^*}+\sum_{u\in N_J(v)} x_u+\sum_{w\in N_S(v)} x_w)-\sum_{v\in V(J)}(k-1)x_v \\
& =|V(J)|+\sum_{v\in V(J)}(d_J(v)-k+1) x_v+\sum_{v\in V(J)} \sum_{w\in N_S(v)} x_w \\
& \leqslant 2k+\eta(V(J))+e(J)+e(J,S) \\
& \leqslant 2k+2k^2+e(J,S). 
\end{align*}
It follows that $e(J,S)>(2k-2)\lambda-(k-1)(2k-2)-2k-2k^2=(2k-2)\lambda-4k^2+2k-2$. Note that $K_k\vee(\frac{m}{k}-\frac{k-1}{2})K_1\in \mathcal{G}(m,F_{2k+2}).$ Hence, $\lambda\geqslant \lambda(K_k\vee(\frac{m}{k}-\frac{k-1}{2})K_1)=\frac{k-1+\sqrt{4m-k^2+1}}{2}>\sqrt{m}\geqslant \frac{3}{2}k^3+2k^2+14k$. Together with Claim~\ref{1^}, we obtain
\begin{align*}
 m & \geqslant d_{R^c}(u^*)+e(R^c)+e(R^c,S)\\
&\geqslant \sum_{J\in \mathcal{J}_5}(|J|+e(J)+e(J,S)) \\
& >(\frac{\lambda}{2k-2}+1)((2k-2)\lambda-4k^2+2k-2+2k^2)\\
&>\lambda^2-(k-1)\lambda-\gamma, 
\end{align*}
 contradicting \eqref{1b}. 
\end{proof}

Denote $\hat{\eta}=\max \{\eta(V(J))\mid J \in \mathcal J_5\}$ for simplicity.
\begin{claim} \label{4^} 
$e(S)\leqslant \hat{\eta}(\frac{\lambda}{2k-2}+1)+\frac{k(k-1)}{2}$.
\end{claim}
\begin{proof}[\bf Proof of Claim~\ref{4^}]
We first show
\[\label{eq:3.006}
 \text{$\eta(V(J))\leqslant 0$\ \  for each\ \  $J\in \mathcal{J}_3.$}
\]
In fact, if $|J|=2k-1$, then $e(J)\leqslant \binom{2k-1}{2}$. Thus by \eqref{7}, we obtain $\eta(V(J))\leqslant(k-1)(2k-1)-(k-1)(2k-1)=0$.
    If $|J|=2k$, then assume without loss of generality that $V(J)=V(C_{2k-1})\cup \{v\}$ and $d_J(v)\geqslant k-1$. One sees that $v$ has just  $k-1$ neighbors in $V(C_{2k-1})$. Otherwise there is a $C_{2k}$ in $J$, contradicting the definition of $\mathcal J_3$. One may also see that $G^*[V(C_{2k-1})]\neq K_{2k-1}$, otherwise combining with $d_J(v)= k-1\geqslant 2$,  there exists a $C_{2k}$ in $J$, a contradiction. Thus $e(J)=e(J-v)+d_J(v)\leqslant (k-1)|J|-1$. Together with \eqref{7}, we have $\eta(V(J))\leqslant e(J)-(k-1)|J|\leqslant -1$ for $J\in\mathcal J_3$ with $|J|=2k$.
    
By Lemmas~\ref{lem3.3a}, \ref{lem3.4a} and~\eqref{eq:3.006}, we have $\eta(V(\hat{J})) \leqslant 0$ for each $\hat{J}\in \mathcal{J}\setminus \mathcal{J}_5$. Thus, $\eta(R^c)=\sum_{J \in \mathcal{J}} \eta(V(J)) \leqslant \sum_{J \in \mathcal J_5}\eta(V(J))\leqslant \hat{\eta}|\mathcal J_5|$. Together with \eqref{6b}, Lemma~\ref{lem3.2a} and Claim~\ref{3^}, we have $e(S) \leqslant \eta(R^c)-\gamma \leqslant \hat{\eta}|\mathcal J_5|+\frac{k(k-1)}{2} \leqslant \hat{\eta}(\frac{\lambda}{2k-2}+1)+\frac{k(k-1)}{2},$ as desired.
 \end{proof}

Now we come back to show Lemma~\ref{lem 5*}.

By Claim~\ref{1^}, one may assume that, for each $J \in\mathcal J_5,$ it is obtained from $K_{2k}$ by deleting $t_J\, (\leqslant k-1)$ edges. So these $t_J$ edges are incident with at most $2t_J$ vertices of $V(J)$. Therefore, there are at least $2 k-2 t_J$ vertices, say  $v_1, v_2, \ldots, v_{2k-2t_J}$, in $V(J)$ such that $d_J(v_1)=\cdots=d_J(v_{2k-2t_J})=2k-1$ and $x_{v_1} \geqslant \cdots \geqslant x_{v_{2k-2t_J}}$.  By Lemma~\ref{lem5.2b}, $N_S(v_1), \ldots, N_S(v_{2k-2t_J})$ and $\bigcup_{i=2k-2t_J+1}^{2k} N_S(v_i)$ are pairwise disjoint. Otherwise, without loss of generality,  suppose that $N_S(v_1)\cup N_S(v_2)\not=\emptyset$, then there exists a copy of $F_{2k+2}$ in $G^*$ with central vertex $v_1$, a contradiction.

Recall that $S=V(G) \setminus N[u^*]$. Let $S_0=\{w \in S| d_S(w)=0\}$ and $S_1=S\setminus S_0$. It is clear that $|S_1| \leqslant 2e(S)$. More precisely, there is no vertex in $S_0$ being the neighbor of $v_i$ for every $i \in\{2, \ldots, 2k-2 t_J\}$. Otherwise, there exists a vertex $w\in N_{S_0}(v_i)$ for some $i\in\{2, \ldots, 2k-2 t_J\},$ and so $N_J(w)=\{v_i\}$. Then 
$G=G^*-v_iw+v_1w$ is an $F_{2k+2}$-free graph with larger spectral radius than $G^*$, a contradiction.  That is to say, $d_S(v_i)=d_{S_1}(v_i)$ for every $i \in\{2, \ldots, 2k-2 t_J\}$. Therefore,  
\begin{align*}
    \lambda\sum_{i=2}^{2k-2t_J}x_{v_i}
    &=\sum_{i=2}^{2k-2t_J}(x_{u^*}+\sum_{u \in N_J(v_i)}x_{u}+\sum_{w\in N_S(v_i)}x_{w})\\
    &\leqslant (2k-2t_J-1)+\sum_{i=2}^{2k-2t_J}d_J(v_i)+\sum_{i=2}^{2k-2t_J}d_{S_1}(v_i)\\
    &\leqslant |S_1|+2k(2k-2t_J-1).
\end{align*}
Recall that $|S_1| \leqslant 2e(S)$. By Claim~\ref{4^}, we obtain
\begin{align*}
\sum_{i=2}^{2 k-2t_J} x_{v_i} &\leqslant \frac{2e(S)+2k(2k-2t_{J}-1)}{\lambda} \\
&\leqslant \frac{2(\hat{\eta}(\frac{\lambda}{2k-2}+1)+\frac{k(k-1)}{2})+2k(2k-2t_{J}-1) }{\lambda}\\
&=\frac{\hat{\eta}}{k-1}+\frac{2\hat{\eta}+5k^2-4kt_J-3k}{\lambda}. 
\end{align*}
Combining with  \eqref{4} gives us
\begin{align*}
 \eta(V(J))
 &=\sum_{u \in V(J)\setminus \{v_2,\ldots,v_{2k-2t_J}\}}(d_J(u)-k+1)x_u+ k\sum_{i=2}^{2k-2t_J}x_{v_i}-e(J) \\
 &\leqslant e(J)-\sum_{u\in \{v_2,\ldots,v_{2k-2t_J}\}}d_J(u)-(k-1)|V(J)\setminus \{v_2,\ldots,v_{2k-2t_J}\}|+k\sum_{i=2}^{2k-2t_J}x_{v_i}\\
& \leqslant \binom{2k}{2}-t_J-(2k-1)(2k-2t_J-1)-(k-1)(2t_J+1)\\
&\ \ \ \ \ +k(\frac{\hat{\eta}}{k-1}+\frac{2\hat{\eta}+5k^2-4kt_J-3k}{\lambda}) \\
& =(2k-1-\frac{4k^2}{\lambda})t_J-2k^2+2k+\frac{k}{k-1}\hat{\eta}+\frac{5 k^3-3 k^2+2 \hat{\eta} k}{\lambda}.
\end{align*}
Bear in mind that $\lambda>\sqrt{m} \geqslant \frac{3}{2}k^3+2k^2+14k$. So we have 
\[ \label{eq:4.2}
\eta(V(J))<(2k-1-\frac{4k}{\frac{3}{2}k^2+2k+14}) t_J-2k^2+2k+(\frac{k}{k-1}+\frac{2}{\frac{3}{2}k^2+2k+14})\hat{\eta}+\frac{5k^2-3k}{\frac{3}{2}k^2+2k+14}.
\]

In order to complete the proof of Lemma~\ref{lem 5*}, it suffices to show the following claim.
\begin{claim}\label{5^}
For each $J\in \mathcal{J}_5$, one has $t_J=k-1$.
\end{claim}
\begin{proof}[\bf Proof of Claim~\ref{5^}]
By Claim~\ref{1^}, $t_J\leqslant k-1$ for each $J\in \mathcal{J}_5$. 
Suppose $0 \leqslant t_J\leqslant k-2$ for some $J \in \mathcal J_5$. By \eqref{eq:3.005}, we have $\hat{\eta}\leqslant k$. Combining with \eqref{eq:4.2} and $k \geqslant 3$, we obtain $\eta(V(J))<\frac{11}{3}+\frac{1}{k-1}+\frac{34k-56}{9k^2+12k+84}-2k<-1$, a contradiction.
\end{proof}

By Claim~\ref{5^}, we have $t_J=k-1$ for each $J\in \mathcal J_5$. Then by~\eqref{7}, we obtain $\eta(V(J))\leqslant e(J)-(k-1)|J| \leqslant 1$ for each $J \in \mathcal J_5$. Thus  $\hat \eta \leqslant 1$. Together with  $k \geqslant 3$ and \eqref{eq:4.2},  $\eta(V(J))<\frac{8}{3}+\frac{1}{k-1}-\frac{2k+44}{9k^2+12k+84}-k<0$, a contradiction. 

This completes the proof.
\end{proof}

Up to now, we know that for all $J\in \mathcal{J}$, one has $\eta(V(J))\leqslant 0$. Combining with~\eqref{6b} and Lemma~\ref{lem3.2a}, we have $e(S)\leqslant \eta(R)-\gamma\leqslant \eta(R^c)-\gamma\leqslant \frac{k(k-1)}{2}$.
Moreover, for every $J\in\mathcal J$, we denote by $\widetilde{J}$ the subgraph of $G^*$ induced by $N_R(V(J))$, where $N_R(V(J))$ is the subset of $R$ in which each vertex has at least one neighbor in $V(J)$. For every member $J\in\mathcal J$, it is clear that $J\subseteq\widetilde{J}$. Consequently, $J$ is the $(k-1)$-core of $\widetilde J$ and $V(J)=(V(\widetilde J))^c$. 

\begin{lem}\label{lem5.3b}
    For each $J\in\mathcal{J}_3$, we have $\eta(V({J}))\leqslant -1$.
\end{lem}
\begin{proof}
Recall that for each $J\in\mathcal J_3$, $J$ contains a longest  cycle of length $2k-1$ and  $2k-1\leqslant|J|\leqslant 2k$. By the proof of \eqref{eq:3.006}, it suffices to show that our result holds for $|J|=2k-1$. 

We first consider $J\in \mathcal J_3$ and $J\ncong K_{2k-1}$. 
For every such $J,$ one sees $e(J)\leqslant e(K_{2k-1})-1= (k-1)|J|-1$. Thus by ~\eqref{7} we have $\eta(V(J))\leqslant e(J)-(k-1)|J|\leqslant -1$, as desired. 

Next we consider $J=K_{2k-1}$. 
If there exists a vertex $v\in V(J)$ such that $d_{G^*}(v)\leqslant k(k+1)$,  then $x_v\leqslant \frac{k(k+1)}{\lambda}$. By \eqref{4}, we have
\begin{align*}
       \eta(V(J))
       &\leqslant \sum\limits_{u\in V(J)\setminus\{v\}}(k-1)x_u+\frac{k^3-k}{\lambda}-e(J)\\
       &\leqslant (2k-2)(k-1)+\frac{k^3}{\lambda}-(2k-1)(k-1)\\
       &= -k+1+\frac{k^3}{\lambda}.
\end{align*}
Together with $\lambda>\sqrt{m}\geqslant \frac{3}{2}k^3+2k^2+14k$ and $k\geqslant 3$, we get $\eta(V(J))<-k+1+\frac{k^2}{\frac{3}{2}k^2+2k+14}<-1$. In what follows, we prove indeed that there exists a vertex $v\in V(J)$ satisfying $d_{G^*}(v)\leqslant k(k+1)$. 

For convenience, let $V(J)=\{v_1,\ldots,v_{2k-1}\}$,  $S'=\bigcup_{i=1}^{2k-1}N_{S}(v_i)$, $S^*_0=\bigcup_{i=1}^{2k-1}N_{S_0}(v_i)$ and $S^*_1=\bigcup_{i=1}^{2k-1}N_{S_1}(v_i)$. Choose two distinct vertices $w_1,w_2$ in $S'$. We are to show the following.
\[\label{eq:3.008}
\text{If $N_J(w_1)\cap N_J(w_2)\neq \emptyset$, then either $N_J(w_1)\subseteq N_J(w_2)$ or $N_J(w_2)\subseteq N_J(w_1)$. }
\] 
In fact, suppose \eqref{eq:3.008} is not true. Then assume $N_J(w_1)\setminus N_J(w_2)\not=\emptyset, N_J(w_2)\setminus N_J(w_1)\not=\emptyset$ and let $v_{1}\in N_J(w_1)\cap N_J(w_2)$. Consequently,  $G^*[\{u^*,w_1,w_2\}\cup V(J)]$ contains an $F_{2k+2}$ with central vertex $v_1$, a contradiction.

By \eqref{eq:3.008}, we can partition $S'$  as $\bigcup_{i=1}^{\ell}S'_i$ such that $\bigcup_{i=1}^{\ell}N_{J}(S'_i)\subseteq V(J)$ and $N_J(S'_i)\cap N_J(S'_j)=\emptyset$ for all $1\leqslant i\neq j\leqslant \ell$. Moreover,  we can assume that $w_i\in S'_i$ and $N_J(w_i)=N_J(S'_i)$ for $1\leqslant i\leqslant \ell$. We proceed by considering the following two possible cases.

{\bf Case 1.} $d_J(w_i)\geqslant 3$ for some $1\leqslant i \leqslant \ell$.

Without loss of generality, we assume $d_J(w_1)\geqslant 3$ and let 
$N_J(w_1)=\{v_1,\ldots,v_d\}$. Then we have $d_J(w)=1$ for any $ w\in S'_1\setminus\{w_1\}$. Otherwise, there is a $w'_i\in S'_1\setminus\{w_1\}$ such that $N_J(w'_i)\subseteq N_J(S'_1)$ and $d_J(w'_i)\geqslant 2$. One may assume that $\{v_1,v_2\}\subseteq N_J(w'_2)$. Then there is an $F_{2k+2}$ in $G^*[\{w_1,w'_2,u^*\}\cup V(J)\}]$, a contradiction. Assume that $x_{v_1}\geqslant \cdots \geqslant x_{v_d}$.
We are to show the following.
\[\label{eq:3.009}
\text{$d_{V(\widetilde J)\setminus V(J)}(v_i)=0$ and $d_{S_0^*\setminus \{w_1\}}(v_i)=0$ for each $2\leqslant i\leqslant d$.}
\]
In fact, if there exists a vertex  $u\in V(\widetilde J)\setminus V(J)$ with $u\sim v_i$ for some $2\leqslant i\leqslant d$. Then there is an  $F_{2k+2}$ with central vertex $v_i$ in  $G^*[\{u^*,u,w_1\}\cup V(J)]$, a contradiction.
If $d_{S_0^*\setminus \{w_1\}}(v_i)> 0$ for some $2\leqslant i\leqslant d$, then $N_J(w)=\{v_i\}$ for each $w\in N_{S_0^*}(v_i)\setminus \{w_1\}$. Let $G=G^*-\{v_iw\mid w\in N_{S_0^*}(v_i)\setminus \{w_1\}\}+\{v_1w\mid w\in N_{S_0^*}(v_i)\setminus \{w_1\}\}$. Then $G$ is an $F_{2k+2}$-free graph and $\lambda(G) > \lambda(G^*)$, a contradiction.

In view of \eqref{eq:3.009}, we obtain  $N_{S_0^*}(v_d)\subseteq \{w_1\}$ and $d_{\widetilde{J}\setminus J}(v_d)=0$.  Thus $d_{G^*}(v_d)\leqslant 1+d_R(v_d)+d_{S'}(v_d) \leqslant 2k+2e(S)\leqslant k(k+1),$ as desired. 

{\bf Case 2.} $d_J(w)\leqslant 2$ for all $w\in S'$.

Recall that $J=K_{2k-1}\in\mathcal J_3$. Let $V(J)=\{v_1,\ldots, v_{2k-1}\}$. Assume that there exist at least two distinct vertices, say $v_i,v_j,$ in $V(J)$ having a common neighbor, say $v,$ in $V(\widetilde{J})\setminus V(J).$ Assume without loss of generality that $x_{v_i}\geqslant x_{v_j}$.
We are to show that $d_{S_0^*}(v_j)=0$. If not, since $G^*$ is $F_{2k+2}$-free, we have $d_{\widetilde{J}}(w)=1$ for each $w\in N_{S_0^*}(v_j)$. Let $G=G^*-\{v_j w\mid w\in N_{S_0^*}(v_j)\}+\{v_iw\mid w\in N_{S_0^*}(v_j)\}$. Clearly, $G$ is $F_{2k+2}$-free and $\lambda(G)>\lambda (G^*)$, a contradiction. Thus, $d_{G^*}(v_j)\leqslant k(k+1)$. 

Now we consider that there is at most one vertex, say $v_{2k-1}$, in $V(J)$ such that it is adjacent to a vertex in $ V(\widetilde{J})\setminus V(J)$. Without loss of generality, we assume that $x_{v_1}=\max \{x_{v_i}\mid i=1,\ldots, 2k-2\}$ and $x_{v_2}+x_{v_3}=\max\{x_{v_i}+x_{v_j}\mid 1\leqslant i\neq j\leqslant 2k-2, N_S(v_i)\cap N_S(v_j)\neq \emptyset\}$. Note that $v_1$ may belong to $\{v_2,v_3\}$. Hence, $d_{S_0^*}(v_i)=0,$ and so $d_{G^*}(v_i)\leqslant k(k+1)$ for each $i\in \{4,\ldots, 2k-2\}$. Otherwise, if there exists some $j\in \{4,\ldots, 2k-2\}$ such that $d_{S_0^*}(v_j)>0$, then, since $G^*$ is $F_{2k+2}$-free, we have $v_{2k-1}\notin N_J(w)$ for all $w\in N_{S_0^*}(v_j)$. Let $G=G^*-\{v w\mid v\in V(J), w\in N_{S_0^*}(v_j),v\sim w\}+\{v_2w,v_3w\mid w\in N_{S_0^*}(v_j),d_J(w)=2\}+\{v_1w\mid w\in N_{S_0^*}(v_j),d_J(w)=1\}$. One sees that $G$ is $F_{2k+2}$-free and $\lambda(G) >\lambda (G^*)$, a contradiction. 

This completes the proof.
\end{proof}
\begin{lem}\label{lem5.4b}
For every member $J\in\mathcal{J}_4$, ${J}$ is a connected component of $G^{*}[R]$ and $\eta(V(J))\leqslant-1$.
\end{lem}
\begin{proof}
According to the definition of $\mathcal J_4$, one sees that, for each $J\in \mathcal J_4$, $J$ contains a cycle $C_{2k}$, and one may label the vertices in $V(J)$ as $v_1,\ldots,v_{2k}$. 
If  there exists a vertex $v\in V(\widetilde J)\setminus V(J)$ such that $N_J(v)\neq\emptyset$,  then there exists a $P_{2k+1}$ in $G^*[V(J)\cup\{v\}]$, a contradiction. Thus $J$ is a connected component of $G^{*}[R]$ for each $J\in\mathcal J_4$. In what follows, we show $\eta(V(J))\leqslant-1$.

Let $S'=\bigcup_{i=1}^{2k}N_{S}(v_i)$, $S^*_0=\bigcup_{i=1}^{2k}N_{S_0}(v_i)$ and $S^*_1=\bigcup_{i=1}^{2k}N_{S_1}(v_i)$. 
Recall that $\eta(V(J))\leqslant e(J)-(k-1)|J|$ and $J$ is obtained from $K_{2k}$ by deleting $t_J$ edges. If $t_J\geqslant k+1$, then $\eta(V(J))\leqslant e(J)-(k-1)|J|\leqslant k(2k-1)-k-1-2k(k-1)=-1$. So we proceed by considering $0\leqslant t_J\leqslant k$. 


{\bf Case 1.} $0\leqslant t_J\leqslant k-1$.\ In this case, there exist at least $2k-2t_J$ vertices, say $v_1,v_2,\ldots,v_{2k-2t_J}$, in $V(J)$, such that $d_J(v_1)=\cdots =d_J(v_{2k-2t_J})=2k-1$ and $x_{v_1}\geqslant\cdots\geqslant x_{v_{2k-2t_J}}$. One sees that $N_S(v_1),\ldots, N_S(v_{2k-2t_J})$ and  $\bigcup_{i=2k-2t_J+1}^{2k}N_S(v_i)$ are pairwise disjoint. 
Otherwise, without loss of generality, assume that there exists a vertex $w_0\in N_S(v_1)\cap \bigcup_{i=2k-2t_J+1}^{2k}N_S(v_i)$. Then by Lemma~\ref{lem5.2b},  $G^*[\{u^*,w_0\}\cup V(J)]$ contains an $F_{2k+2}$ with central vertex $v_1$, a contradiction. We claim $N_{S_0}(v_i)=\emptyset$, i.e., $d_S(v_i)=d_{S_1}(v_i)$ for each $i\in \{2,\ldots, 2k-2t_J\}$. Otherwise, there exists a vertex $w\in N_{S_0}(v_i)$ for some $2\leqslant i \leqslant 2k-2t_J$. By~Lemma~\ref{lem5.1}, $G^*-v_i w+v_1w$ is $F_{2k+2}$-free and has a larger spectral radius than $G^*$, a contradiction.   
Recall that $|S_1|\leqslant 2e(S)$. Thus,  
\begin{align}
\lambda\sum_{i=2}^{2k-2t_J}x_{v_i}
&=\sum_{i=2}^{2k-2t_J}(x_{u^*}+\sum_{u\in N_J(v_i)}x_u+\sum_{w\in N_{S}(v_i)}x_w)\notag\\
&\leqslant(2k-2t_J-1)+\sum_{i=2}^{2k-2t_J}d_J(v_i)+\sum_{i=2}^{2k-2t_J}d_S(v_i)\notag\\
&\leqslant (2k-2t_J-1)+(2k-1)(2k-2t_J-1)+|S_1|\notag\\
&\leqslant 2e(S)+2k(2k-2t_J-1)\notag\\
&\leqslant k(k-1)+2k(2k-2t_J-1). \tag{As $e(S)\leqslant \frac{k(k-1)}{2}$}
\end{align}
Hence,
$$
\sum_{i=2}^{2k-2t_J}x_{v_i}\leqslant \frac{k(k-1)+2k(2k-2t_J-1)}{\lambda}.
$$

Recall that $J$ is obtained from $K_{2k}$ by deleting $t_J$ edges 
and $d_J(v_2)=\cdots =d_J(v_{2k-2t_J})=2k-1$. 
Hence, we have  
$
\sum_{u\in V(J)\setminus\{v_2,\ldots,v_{2k-2t_J}\}}(d_J(u)-k+1)x_u\leqslant 2e(J)-\sum_{u\in\{v_2,\ldots,v_{2k-2t_J}\}}d_J(u)-(k-1)(2t_J+1)\leqslant 2k(2k-1)-2t_J-(2k-1)(2k-2t_J-1)-(k-1)(2t_J+1)=(2k-2)t_J+k$.
By \eqref{4}, we have
\begin{align*}
       \eta(V(J))
       &\leqslant \sum_{u\in V(J)\setminus\{v_2,\ldots,v_{2k-2t_J}\}}(d_J(u)-k+1)x_u+\frac{k^2(k-1)+2k^2(2k-2t_J-1)}{\lambda}-e(J)\\
       &\leqslant (2k-2)t_J+k-k(2k-1)+t_J+\frac{5k^3-4k^2t_J-3k^2}{\lambda}\\
       &= (2k-1-\frac{4k^2}{\lambda})t_J-2k^2+2k+\frac{5k^3-3k^2}{\lambda}.  
  \end{align*}
Recall that $ 0\leqslant t_J\leqslant k-1$,  $\lambda>\sqrt{m}\geqslant \frac{3}{2}k^3+2k^2+14k$ and $k\geqslant 3.$ Consequently, 
$$
\eta(V(J)) 
\leqslant (2k-1-\frac{4k}{\frac{3}{2}k^2+2k+14})(k-1)-2k^2+2k+\frac{5k^2-3k}{\frac{3}{2}k^2+2k+14}
<-1.
$$ 

{\bf Case 2.} $t_J=k$.\ In this case, by \eqref{7}, we have $\eta(V(J))\leqslant e(J)-(k-1)|J|=k(2k-1)-k-(k-1)2k=0$. 

If  there exists a vertex $v_i\in V(J)$ satisfying $d_J(v_i)\geqslant 2k-2$ and $d_{G^*}(v_{i})\leqslant  k(k+1)$,
we obtain  $x_{v_{i}}\leqslant\frac{k(k+1)}{\lambda}$. 
Consequently,  $\sum_{u\in V(J)\setminus \{v_{i}\}}(d_J(u)-k+1)x_u\leqslant \sum_{u\in V(J)\setminus \{v_{i}\}}d_J(u)-(k-1)(2k-1)\leqslant 2e(J)-d_J(v_{i})-(k-1)(2k-1)\leqslant 2k^2-3k+1$.
By \eqref{4}, we have
 \begin{align*}
       \eta(V(J))
       &\leqslant \sum\limits_{u\in V(J)\setminus\{v_{i}\}}(d_J(u)-k+1)x_u+\frac{k^3+k^2}{\lambda}-e(J)\\
       &\leqslant 2k^2-3k+1+\frac{k^3+k^2}{\lambda}-k(2k-1)+k\\
       &= -k+1+\frac{k^3+k^2}{\lambda}.
  \end{align*}
Recall that $\lambda>\sqrt{m}\geqslant \frac{3}{2}k^3+2k^2+14k$. Then $\eta(V(J))<-k+1+\frac{k^2+k}{\frac{3}{2}k^2+2k+14}<-1.$

Next, we  show that there indeed exists a vertex $v_i\in V(J)$ such that $d_J(v_i)\geqslant 2k-2$ and  $d_{G^*}(v_i)\leqslant k(k+1)$ by considering the following two  subcases.

{\bf Subcase 2.1.} There exist at least two distinct vertices $v_i, v_j \in V(J)$ such that $d_J(v_i)=d_J(v_j)= 2k-1$. 
In this subcase, assume without loss of generality that $d_J(v_1)=d_J(v_2)=2k-1$ where $x_{v_1}\geqslant x_{v_2}$. By a similar discussion as that in Case~1, one sees that $N_S(v_1), N_S(v_2)$ and $\bigcup_{i=3}^{2k}N_S(v_i)$ are pairwise disjoint.  
If  $d_{S_0}(v_2)\neq 0$ and $w\in N_{S_0}(v_2)$, then let $G=G^*-v_2w+v_1w$. Clearly $G$ is $F_{2k+2}$-free and has larger spectral radius than $G^*$, a contradiction.
Thus $d_{S_0}(v_2)=0$. 
Recall that $e(S)\leqslant\frac{k(k-1)}{2}$. Then  $d_{G^*}(v_2)\leqslant 1+2k-1+2e(S)\leqslant k(k+1)$. 

{\bf Subcase 2.2.} There exists at most one vertex $v \in V(J)$ with $d_J(v)=2k-1$. In this subcase, we may partition $V(J)$ as $V_1\cup V_2$, where $|V_1|=|V_2|=k$, and $J[V_1] \cong J[V_2] \cong K_k$. We proceed by showing the following claims to complete our proof.
\begin{claim}\label{3*}
All the vertices $w \in N_{S_0}(V(J))$ with $d_{J}(w)=1$ have a common neighbor in $V(J)$.
\end{claim}
\begin{proof}[\bf Proof of Claim~\ref{3*}]
Let $v$ be in $V(J)$ with $x_v=\max _{u \in V(J)} x_u$, and let $w$ be a vertex in $N_{S_0}(V(J))$ with $d_J(w)=1$. If $w \nsim v$, then let $G=G^*-w v'+w v$, where $v'$ is the unique neighbor of $w$ in $V(J)$. Clearly $G$ is $F_{2k+2}$-free and has larger spectral radius than $G^*$, a contradiction.
\end{proof}

\begin{claim}\label{4*}
All the vertices $w \in N_{S_0}(V(J))$ with $d_J(w)=2$ have a common neighborhood in $V(J)$.
\end{claim} 
\begin{proof}[\bf Proof of Claim~\ref{4*}]
Suppose there are two distinct vertices $w$ and $w'$ in $N_{S_0}(V(J))$ with $d_{J}(w)=d_{J}(w')=2$ satisfying $N_J(w) \neq N_{J}(w')$. Assume $\sum_{v \in N_{J}(w)} x_v \geqslant \sum_{v \in N_{J}(w')} x_v$. Then let $G=G^*-\{w'v\mid v\in N_{J}(w')\}+\{w'v\mid v\in N_{J}(w)\}$. Clearly, $G$ is $F_{2k+2}$-free and has a larger spectral radius, a contradiction.
\end{proof}

By Claims~\ref{3*} and~\ref{4*}, we may assume all vertices in $N_{S_0}(V(J))$ with only one neighbor in $V(J)$ (if there exists) have common neighbor $v_{i_1}\in V(J)$, and all vertices in  $N_{S_0}(V(J))$ with exactly two neighbors in $V(J)$ (if there exist) have common neighborhood $\{v_{i_2}, v_{i_3}\}$. Note that $v_1$ may belong to $\{v_2,v_3\}.$ Then for each $v_i\in V(J) \setminus \{v_{i_1}, v_{i_2}, v_{i_3}\}$ and each $w\in N_{S_0}(v_i)$, we obtain  $d_J(w) \geqslant 3$. Further on we have the following claim.
\begin{claim}\label{5*}
Let $v,v'$ be two distinct vertices of $J$ satisfying $v\not\sim v'$ and $d_J(v)=2k-2$. Then, for any distict vertices $v_i, v_j \in V(J)\setminus \{v, v'\}$, $G^*[\{u^*\}\cup (V(J) \setminus \{v, v'\})]$ contains a $v_i v_j$-path of length $2k-2$.
\end{claim}
\begin{proof}[\bf Proof of Claim~\ref{5*}]
Recall that $V_1\cup V_2$ is a partition of $V(J)$. Hence one may assume, without loss of generality, that $v\in V_1$ and $v'\in V_2$. By the symmetry of $v_i$ and $v_j$, we proceed by considering the following two cases. 

$\bullet$ $v_i\in V_1, v_j\in V_2$. Then $v_i v_{i_1'}\cdots, v_{i_{k-2}'}u^*v_{j_{k-2}}\cdots v_{j_1} v_j$ is a desired $v_iv_j$-path, where $v_{i_1'}, \ldots, v_{i_{k-2}'}$ $\in V_1\setminus\{v, v_i\}$ and $v_{j_1}, \ldots, v_{j_{k-2}}\in V_2\setminus\{v', v_j\}$.

$\bullet$  $v_i, v_j \in V_1$ or $v_i, v_j \in V_2$. Here we only consider the former. Choose some vertex $v_{j_{k-1}}\in V_2\setminus\{v'\}$ such that $v_{j_{k-1}}\sim v_j$. Then $v_iv_{i_1}\cdots v_{i_{k-3}} u^* v_{j_1} \cdots v_{j_{k-1}} v_j$ is a desired $v_iv_j$-path, where $v_{i_1},\ldots, v_{i_{k-3}}\in V_1 \setminus\{v, v_i, v_j\}$ and $v_{j_1}, \ldots, v_{j_{k-1}}\in V_2\setminus \{v'\}$. In particular, $v_i=v_{i_0}$ if $k=3$. 
\end{proof}

Next we come back to continue the proof for Subcase 2.2.

If there is a vertex $w \in N_{S_0}(V(J))$ with $d_J(w)\geqslant 3$ and there is a vertex $v_i\in N_J(w)\setminus\{v_{i_1}, v_{i_2}, v_{i_3}\}$  such that $d_J(v_i)\geqslant 2k-2$, and $v_i$ is adjacent to at least two vertices in $N_J(w)$, then $w$ is the unique vertex in $S_0$ satisfying $w\sim v_i$. Otherwise, suppose $w'\in N_{S_0}(v_i)\setminus\{w\}$, by Claims~\ref{3*} and \ref{4*},  we have $d_J(w')\geqslant 3$, and so by Claim~\ref{5*},  $G^*[\{u^*,w, w'\}\cup N_J[v_i]]$ contains an $F_{2k+2}$ with central vertex $v_i$, a contradiction. Therefore, $N(v_i)\subseteq
\{u^*,w\}\cup V(J)\cup S_1$, and so $d_{G^*}(v_i)\leqslant k(k+1)$.

Now we consider $d_J(w)\geqslant 3$ for all $w$ in $N_{S_0}(V(J))$ and each vertex $v_i\in N_J(w)\setminus\{v_{i_1},v_{i_2},v_{i_3}\}$ with $d_J(v_i)\geqslant 2k-2$ is adjacent to at most one vertex in $N_J(w)$. For such vertex $v_i$, if $d_{S_0}(v_i)\leqslant 1$, then $d_{G^*}(v_i)\leqslant k(k+1)$. If $d_{S_0}(v_i)\geqslant 2$, then all vertices in $N_{S_0}(v_i)$ have a common neighborhood of size\footnote{Here we use \textit{size} to denote the cardinality of a set} three including $v_i$ and the unique vertex, say $v_{i^*}$, in $V(J)$ satisfying $v_{i^*}\nsim v_i$. Otherwise, by Claim~\ref{5*}, there is an $F_{2k+2}$ in $G^*$ with central  vertex $v_i$, a contradiction. Take $w, w'\in N_{S_0}(v_i)$ and let $N_J(w)=N_J(w')=\{v_i, v_{i^*}, v_{i'}\}$. If $v_{i^*}\sim v_{i'}$ and $d_J(v_i')\geqslant 2k-2$, then by Claim~\ref{5*}, there exists an $F_{2k+2}$ in $G^*$ with central  vertex $v_{i'}$, a contradiction. If $v_{i^*}\nsim v_{i'}$ or $d_J(v_{i'})=2k-3$, then $N_J(w)$ contains a vertex of degree $2k-3$ in $J$. In the following, we consider that for all $w \in N_{S_0}(V(J))$ with $d_J(w) \geqslant 3$, $N_J(w)$ contains a vertex of degree less than $2k-2$ in $J$. Then there is a vertex, say $v_{i_4}$,  of degree $2k-1$ in $J$, a unique vertex, say $v_{i_5}$, of degree less than $2k-2$ in $J$ satisfying $d_J(v_{i_5})=2k-3$. If $N_{S_0}(v_{i_4})=\emptyset$, then $d_{G^*}(v_{i_4})\leqslant 1+2k-1+|S_1|\leqslant k(k+1)$. If $N_{S_0}(v_{i_4})\neq \emptyset$, then $v_{i_4}=v_{i_1}$. Otherwise, by Claim~\ref{5*}, $G^*$ contains an $F_{2k+2}$ with central vertex $v_{i_4}$, a contradiction. 

Furthermore, we have the following two claims.
\begin{claim}\label{6*}
  Let $w\in N_{S_0}(V(J))$ satisfying $d_J(w)\geqslant 3$. Then either $d_J(w)=3$ or $d_{G^*}(v)\leqslant k(k+1)$ for some $v \in V(J)$.
\end{claim}
\begin{proof}[\bf Proof of Claim~\ref{6*}]
If there is a vertex $w\in N_{S_0}(V(J))$ such that $d_J(w)\geqslant 4$, then let $v_i\in N_J(w)\setminus\{v_{i_2}, v_{i_3}, v_{i_5}\}$. One sees that $v_i$ has at least two neighbors in $N_J(w)$. If $d_{S_0}(v_i)\geqslant 2$, by Claim~\ref{5*} there is an $F_{2k+2}$ in $G^*$ with central  vertex $v_i$, a contradiction. If $d_{S_0}(v_i)=1$, then $d_{G^*}(v_i)\leqslant 1+2k-1+|S_1|\leqslant k(k+1)$, as desired.
\end{proof}
\begin{claim}\label{7*}
 Let $w, w'$ be two distinct vertices in $N_{S_0}(V(J))$ satisfying $d_J(w)=d_J(w')=3$. Then either $N_J(w)=N_J(w')$ or $N_J(w)\cup N_J(w')=\{v_{i_5}\}$. Furthermore, if $\{v_{i_2},v_{i_3}\}\nsubseteq N_J(w)$ and $\{v_{i_2},v_{i_3}\}\nsubseteq N_J(w')$,  then either  $N_J(w)=N_J(w')$ or $d_{G^*}(v)\leqslant k(k+1)$ for some $v \in V(J)$.
\end{claim}
\begin{proof}[\bf Proof of Claim~\ref{7*}]
Suppose $N_J(w)\neq N_J(w')$. Then we have $v_{i_5}\in N_J(w)\cap N_J(w')$. Assume that $v_i$ is a common vertex of $N_J(w)$ and $N_J(w')$ other than $v_{i_5}$. Then $G^*$ contains an $F_{2k+2}$ with central vertex $v_i$, a contradiction. Hence $N_J(w) \cap N_J(w')=\{v_{i_5}\}$. 

Now we  consider the second part of this claim. Suppose $N_J(w) \neq N_J(w')$, then we may assume  $N_J(w)=\{v_{i_5}, v_{j_1}, v_{j_2}\}$ and $N_J(w')=\{v_{i_5}, v_{j_3}, v_{j_4}\}$ with $\{v_{j_1},v_{j_2}\}\cap \{v_{j_3},v_{j_4}\}=\emptyset$.  We first consider $\min\{d_{S_0}(v_{j_i})|i=1,2,3,4\}\geqslant 2$. Assume without loss of generality that $x_{v_{j_1}}+x_{v_{j_2}}\geqslant x_{v_{j_3}}+x_{v_{j_4}}$. Let $G=G^*-\{w'v_{j_3}, w'v_{j_4}\}+\{w'v_{j_1}, w'v_{j_2}\}$. Then $\lambda(G)>\lambda(G^*)$. 
On the other hand, one may see that $G$ is $F_{2k+2}$-free, which derives a contraction to the choice of $G^*$. In fact, if $G$ contains an $F_{2k+2}$, then this $F_{2k+2}$ must contain $w'$. Note that $\min\{d_{S_0}(v_{j_i})|i=1,2,3,4\}\geqslant 2$ and $\{v_{j_1},v_{j_2}\}\nsubseteq\{v_{i_1},v_{i_2},v_{i_3}\}$. 
Without loss of generality, we assume $v_{j_1}\notin\{v_{i_1},v_{i_2},v_{i_3}\}$. Then there is a vertex $w''\in N_{S_0}(V(J))\setminus\{w\}$ such that $d_J(w'')=3$ and $v_{j_1}\in N_J(w'')$. Hence, 
$\{v_{i_5},v_{j_1}\}\subseteq N_J(w'')$. Now 
$\{v_{i_5},v_{j_1}\}\subseteq N_J(w)\cap N_J(w')$, and so 
$N_J(w'')= N_J(w)=\{v_{i_5},v_{j_1},v_{j_2}\}$. Hence, in $G$, one has $N_J(w)=N_J(w'')=N_J(w')=\{v_{i_5},v_{j_1},v_{j_2}\}$. Since an $F_{2k+2}$ in $G$ may not contain $w,w'$  and $w''$ simultaneously. After replacing $w'$ with $w$ or $w'$, we find another $F_{2k+2}$ in $G$, which is also in $G^*$. That is, $G^*$ contains an $F_{2k+2}$, a contradiction. Therefore, $G$ is $F_{2k+2}$-free. The remaining case is  $\min\{d_{S_0}(v_{j_i})|i=1,2,3,4\}\leqslant 1$. Assume  $d_{S_0}(v_{j_1})=\min\{d_{S_0}(v_{j_i})|i=1,2,3,4\}$. Then 
$d_{G^*}(v_{j_1}) \leqslant k(k+1)$, as desired.
\end{proof}

Now we come back to complete the proof for Subcase 2.2. 

By Claims~\ref{3*},~\ref{4*},~\ref{6*} and~\ref{7*}, if $k\geqslant 4$, then there is a vertex $v\in V(J)$ such that $d_J(v)=2k-2$ and  $d_{S_0}(v)\leqslant 1$, and so $d_{G^*}(v)\leqslant k(k+1)$. We proceed by considering $k=3$. In this subcase, if $N_{S_0}(v_{i_1})=\emptyset$, then $d_{G^*}(v)\leqslant k(k+1)$. If $N_{S_0}(v_{i_1})\neq\emptyset$, then by Claims~\ref{3*} and~\ref{5*}, we know that $N_{S_0}(v_{i_1})=\{w\in S_0|d_J(w)=1\}$. Hence $x_{v_{i_1}}=\max\{x_v|v\in V(J)\}$, otherwise, suppose there is a vertex $v_{i_j}\in V(J)$ such that $x_{v_{i_j}}>x_{v_{i_1}}$. Construct $G=G^*-\{v_{i_1}w|w\in N_{S_0}(v_{i_1})\}+\{v_{i_j}w|w\in N_{S_0}(v_{i_1})\}$, then $G$ is $F_{2k+2}$-free and has a larger spectral radius than $G^*$, a contradiction. 
Recall that $\eta(V(J))=\sum_{u\in V(J)}(d_J(u)-2)x_u-e(J)$. If $x_{v_{i_1}}\leqslant \frac{11}{12}$, then 
\begin{align*}
\eta(V(J)) &\leqslant \sum_{u \in V(J)} (d_{J}(u)-2)x_{v_{i_1}} - e(J) \\
&= (2e(J)-12)x_{v_{i_1}} - e(J) \\
&= 12x_{v_i} - 12 \\
&\leqslant -1,
\end{align*}
as desired. If $\sum_{u \in V(J)\setminus \{v_{i_1}\}} x_{u} \leqslant 4 $, then
\begin{align*}
\eta(V(J)) &\leqslant \sum_{u \in V(J)\setminus\{v_{i_1}\}} x_{u} + x_{v_{i_1}} + \sum_{u \in V(J)} (d_{J}(u)-3)x_{u} - e(J) \\
&\leqslant 4 + 1 + (2e(J)-18) - e(J) \\
&= -1,
\end{align*}
as desired. 

In the following, we assume $x_{v_{i_1}} > \frac{11}{12}$ and $\sum_{u \in V(J)\setminus\{v_{i_1}\}} x_{u} > 4$.
By $\lambda \mathbf{x} = A(G^{*})\mathbf{x}$, one has
\begin{align*}
\lambda x_{v_{i_1}} &= \sum_{v \in N(v_{i_1})} x_{v} \\
&= x_{u^*} + \sum_{u \in N_J(v_{i_1})} x_{u} + \sum_{w \in N_{S_0}(v_{i_1})} x_{w}+\sum_{w'\in N_{S_1}(v_{i_1})}x_{w'} \\
&\leqslant 1 + 5x_{v_i} + \sum_{w \in N_{S_0}(v_{i_1})} x_{w} + |S_1|.
\end{align*}
Then we obtain 
\begin{align*}
\sum_{w \in N_{S_0}(v_{i_1})} x_{w} &\geqslant (\lambda-5)x_{v_{i_1}} - 1-|S_1|
\geqslant (\lambda - 5)x_{v_{i_1}} -1- 6 
> \frac{11}{12}\lambda - \frac{139}{12}.
\end{align*}

On the other hand, for each $ w \in N_{S_0}(v_{i_1}) $, since \( d_{J}(w) = 1 \), one has
\begin{align*}
\lambda x_{w} = \sum_{u \in N(w)} x_{u}
=\sum_{u \in N_R(w)} x_{u}
\leqslant \sum_{u \in R} x_{u} -\sum_{u \in V(J)\setminus \{v_{i_1}\}} x_{u}<\lambda - 4,
\end{align*}
and so
\begin{align*}
\sum_{w \in N_{S_0}(v_{i_1})} x_{w} < |N_{S_{0}}(v_{i_1})| - \frac{4}{\lambda}|N_{S_0}(v_{i_1})|.
\end{align*}
This leads to
\begin{align*}
|N_{S_{0}}(v_{i_1})| > \frac{1}{1 - \frac{4}{\lambda}} \sum_{w \in N_{S_0} (v_{i_1})} x_{w} 
\geqslant \frac{1}{1 - \frac{4}{\lambda}} ( \frac{11}{12}\lambda - \frac{139}{12})
=\frac{\lambda}{\lambda - 4} ( \frac{11}{12}\lambda - \frac{139}{12}).
\end{align*}
Then we have 
\begin{align*}
\sum_{w \in N_{S_0}(v_{i_1})} x_{w} 
&< |N_{S_{0}}(v_{i_1})| - \frac{4}{\lambda}\frac{\lambda}{\lambda - 4} ( \frac{11}{12}\lambda - \frac{139}{12})\\
&=|N_{S_0}(v_{i_1})| - \frac{4}{\lambda - 4} ( \frac{11}{12}\lambda - \frac{139}{12})\\
&=|N_{S_0}(v_{i_1})|-\frac{11}{3}+\frac{95}{3(\lambda-4)}.
\end{align*}
Note that $\lambda > \sqrt{m} \geqslant 99$, so we have 
\begin{align*}
\sum_{w \in N_{S_0}(v_{i_1})} x_{w} &< |N_{S_{0}}(v_{i_1})| - \frac{11}{3} + \frac{1}{3}=|N_{S_0}(v_{i_1})|-\frac{10}{3}.
\end{align*}
Now by \eqref{5.0b}, we have 
\begin{align*}
\lambda^2 - 2\lambda &= d_{G^{*}}(u^*) + \eta(R) + e(R) + \sum_{w \in S} d_{R}(w)x_{w} \\
&= d_{G^*}(u^*) + \eta(R) + e(R) +\sum_{w\in N_{S_0}(v_{i_1})}x_w +\sum_{w\in N_{S_0}(v_{i_1})}(d_R(w)-1)x_w+\sum_{w\in S\setminus N_{S_0}(v_{i_1})}x_w\\
&\leqslant d_{G^*}(u^*) + \eta(R) + e(R) +e(R,S)-|N_{S_0}(v_{i_1})|+\sum_{w\in N_{S_0}(v_{i_1})}x_w\\
&<d_{G^*}(u^*)+\eta(R)+e(R)+e(R, S) - \frac{10}{3}\\
&=\eta(R)+m-e(S)-\frac{10}{3}.
\end{align*}
Combining with \eqref{1b}, we get $\eta(R) > \frac{1}{3} + e(S) > 0$, a contradiction.

Therefore, we deduce that $J$ contains a vertex $v$ satisfying $d_{G^*}(v)\leqslant k(k+1)$ for $k\geqslant 3$.  Consequently, $\eta(V(J))\leqslant -1$ for each $J\in \mathcal{J}_4$, as desired.
\end{proof}
\subsection{\normalsize On the cardinalities of $\mathcal{J}_1, \mathcal{J}_2, \mathcal{J}_3$ and $\mathcal{J}_4$ }
In this subsection, we determine the cardinalities, respectively, for $\mathcal J_1, \mathcal J_2, \mathcal J_3$ and $\mathcal J_4$. For convenience,  we give a partition of $R\setminus R^{c}$: Let $Q = \{u \in R\setminus R^{c}| d_R(u)\leqslant k-2\}$ and $P = R \setminus (R^{c}\cup Q)$. Denote $|P|=p$, $|Q|=q$ for simplicity. 
\begin{lem}\label{5.5b}
 $|\mathcal{J}_1|=1$ and $|\mathcal{J}_2|=|\mathcal{J}_3|=|\mathcal{J}_4|=0$.
\end{lem}
\begin{proof}
If $|\mathcal J_1|\geqslant 2$,  by Lemma~\ref{lem3.3a},  we obtain $\eta(R)\leqslant (\gamma+1)|\mathcal J_1|<\gamma$, contradicting \eqref{6b}. 
Then we suppose  that $|\mathcal{J}_1|=0$. Then $\mathcal{J}= \mathcal{J}_2\cup \mathcal{J}_3\cup \mathcal{J}_4$ and so $|J|\leqslant 2k$ for each $J\in \mathcal{J}$. Further by Lemmas~\ref{lem3.4a}, \ref{lem5.3b} and \ref{lem5.4b}, we know that  $\eta(V(J))\leqslant -1$ for each $J\in \mathcal{J}$. Therefore, $|\mathcal{J}|\leqslant \frac{k(k-1)}{2}$.  Otherwise, by Lemma~\ref{lem3.2a},  $\eta(R) \leqslant \eta(R^c) <-\frac{k(k-1)}{2}=\gamma $, contradicting \eqref{6b}. Combining with $R^c=\cup_{J\in \mathcal{J}}V(J)$, one has $|R^c|\leqslant k^2(k-1)$ and $e(R^c)\leqslant (k-\frac{1}{2})|R^c|$.  
    
 By the definition of $(k-1)$-core, $P\cup Q$ admits a vertex ordering $u_1,\ldots,u_{p+q}$ such that $d_{R_{i}}(u_{i})\leqslant k-2$ for $i\in\{1,\ldots,p+q\}$, where $R_1 = R$ and $R_i = R_{i-1}\setminus \{u_{i-1}\}$ for $i\geqslant2$. Since the $k$-core is well-defined and $d_{R}(u)\leqslant k-2$ for each $u \in Q$, we may assume that $Q=\{u_{1},\ldots,u_{q}\}$ and $P=\{u_{q+1},\ldots,u_{p+q}\}$. By the definition of $(k-1)$-core, one may easily see $e(P)+e(P, R^c)=\sum_{i=q+1}^{p+q}d_{R_i}(u_i)$, and so $e(P)+e(P, R^c)\leqslant (k-2)p$.

Observe $\sum_{u\in P\cup R^c}d_R(u)\leqslant e(P\cup R^c)+e(R)$, where $e(P\cup R^c)=e(P)+e(P, R^c)+e(R^c)$. Then
\begin{align*}
\eta(R)& =\sum_{u\in Q}(d_R(u)-k+1)x_u+\sum_{u\in P\cup R^c}(d_R(u)-k+1)x_u-e(R)  \\
&\leqslant-\sum_{u\in Q}x_u+e(P)+e(P, R^c)+e(R^c)-(k-1)(p+|R^c|).
\end{align*}
It follows that $\eta(R)\leqslant-\sum_{u\in Q}x_u-p+\frac{|R^c|}{2}$.  In view of \eqref{6b}, we know that $\eta(R)\geqslant \gamma=-\frac{k(k-1)}{2}$. Thus $\sum_{u\in Q}x_{u}\leqslant\frac{k^{2}-k}{2}-p+\frac{|R^c|}{2}$. Recall that $|R^c|\leqslant k^2(k-1).$ Then
\begin{align*}
        \lambda=\lambda x_{u^*}=\sum_{u\in Q}x_u+\sum_{u\in P}x_u+\sum_{u\in R^c}x_u\leqslant\frac{3k^3-2k^2-k}{2},
\end{align*}
which contradicts $\lambda> \sqrt{m}\geqslant \frac{3}{2}k^3+2k^2+14k $. Therefore, $|\mathcal{J}_1|=1$, as desired. 

We now prove $|\mathcal{J}_2|=0$. By Lemmas~\ref{lem3.3a},~\ref{lem3.4a},~\ref{lem5.3b} and~\ref{lem5.4b},  one has  $\eta(V( J))\leqslant -(k-1)$ for each $J\in\mathcal J_2$,  $\eta(V(J))\leqslant \gamma+1$ for each $J\in \mathcal{J}_1$  and $\eta(V( J))\leqslant -1$ for each $J\in\mathcal J_3\cup \mathcal J_4$. 
As $|\mathcal{J}_1|\geqslant 1$, if  $|\mathcal J_2|\geqslant 1$, then  $\eta(R)\leqslant \eta(R^c)\leqslant \sum_{J\in\mathcal J_1}\eta(V(J))+\sum_{ J\in\mathcal J_2}\eta(V(J))<\gamma$, contradicting \eqref{6b}. Thus we obtain $|\mathcal J_2|=0$. 
   
Now we show that $|\mathcal J_3|=|\mathcal J_4|=0$.  If $|\mathcal{J}_3|+|\mathcal{J}_4|\geqslant 2$, then  combining 
$|\mathcal{J}_1|\geqslant 1$ with Lemmas~\ref{lem3.3a},~\ref{lem5.3b} and~\ref{lem5.4b}, we have  $\eta(R)\leqslant \eta(R^c)\leqslant \gamma-1<\gamma$, a contradiction.      
If $|\mathcal{J}_3|+|\mathcal{J}_4|=1$,  we may assume, without loss of generality, that $|\mathcal J_3|=1, |\mathcal J_4|=0$. Let $J$ be the unique element in $\mathcal{J}_1$.  

If $J\ncong  S_{|J|,k-1}^+$, by Lemma~\ref{lem3.3a} we obtain $\eta(V(J))\leqslant \gamma$. Thus $\eta(R)\leqslant\eta(R^c)\leqslant \eta(V(J))+\sum_{\hat J\in\mathcal J_3}\eta(V(\hat J))\leqslant \gamma-1<\gamma$, a contradiction. 

If $J\cong S_{|J|,k-1}^+$ and there exists a vertex $v\in V(J)$ with $d_J(v)\geqslant k$ and  $x_v<1$, then by \eqref{4}, we have $\eta(V(J))<\sum_{u\in V(J)}(d_J(u)-k+1)-e(J)=\gamma+1$.  By Lemmas~\ref{lem3.2a} and~\ref{lem5.3b}, we obtain $\eta(R)\leqslant\eta(R^c)\leqslant \eta(V(J))+\sum_{\hat J\in\mathcal J_3}\eta(V(\hat J))< (\gamma+1)-1=\gamma$,  a contradiction.

If $J\cong S_{|J|,k-1}^+$ and $x_v=1$ for each  $v\in V(J)$ with $d_J(v)\geqslant k$,  by~\eqref{4}  one sees $\eta(V(J))=\gamma+1$. By Lemma~\ref{lem5.3b}, we obtain $\eta(R^c)\leqslant \eta(V(J))+\sum_{\hat J\in\mathcal J_3}\eta(V(\hat J))\leqslant \gamma$. Recall that $\eta(R^c)\geqslant \eta(R)\geqslant \gamma $. Then  $\eta(R^c)=\eta(R)=\gamma$.  Hence, by Lemma~\ref{lem3.2a}  one has $R=R^c$,   and  by \eqref{5b} and \eqref{6b}, we also have $e(S)=0$ and $x_w=1$ for each $w\in S$. In what follows, we show that $S=\emptyset.$ 

Suppose to the contrary that there exists a vertex $w\in S$ such that $x_w=x_{u^*}=1$. Combining $e(S)=0$,  we have $N(w)=N(u^*)$. Thus $G^*[\{u^*,w\}\cup R]$ contains a $P_{2k+1}=wv_1v_2u_2v_3\ldots u_{k-1}v_k\linebreak u^*v_{k+1}$, where $v_1\sim v_2$ and each vertex in $\{u_2,\ldots,u_{k-1}\}$ is of degree $|J|-1$ in $J$. Then $F_{2k+2}$ is a subgraph of  $G^*[\{u^*,w\}\cup R]$, a contradiction. So we obtain $S=\emptyset$. 
Note  that there is a vertex $v\in V(J)$  such that $d_J(v)=k$. 
Thus $d(v)=1+k$. Together with $\lambda>\sqrt{m}\geqslant \frac{3}{2}k^3+2k^2+14k$ and $k\geqslant 3$,   we have $x_{v}\leqslant \frac{k+1}{\lambda}<1$, a contradiction.  
Thus $|\mathcal J_3|=|\mathcal J_4|=0$, as desired. 
     
This completes the proof.
\end{proof}

\section{\normalsize Proofs of Theorem \ref{thm2.1}, Corollaries~\ref{thm1.3} and~\ref{thm2.2}}\label{sec-4}

In this section, we give the proofs of Theorems~\ref{thm2.1}, Corollaries~\ref{thm1.3} and~\ref{thm2.2}, respectively. Theorem~\ref{thm2.1} determines the unique graph among $\mathcal G(m,F_{2k+2})$ having the largest spectral radius, which deduces Conjecture~\ref{conj-6} directly. Corollary~\ref{thm1.3} characterizes the unique graph among $\mathcal G(m,F_{k,3})$ having the largest spectral radius, which deduces Conjecture~\ref{cj1-5} directly. Corollary~\ref{thm2.2} identifies the graphs among $\mathcal{G}(m,\theta_{1,p,q})$ for $q \geqslant p\geqslant 3$ having the largest spectral radius, which resolves Problem~\ref{pb-2} for $q+p\geqslant 7$.  

\begin{proof}[\bf Proof of Theorem~\ref{thm2.1}]

Recall $\mathcal J_1=\{J\in\mathcal J:|J|\geqslant2k+1\}$. In view of Lemma~\ref{5.5b}, $\mathcal J_1=\{R^c\}$ and then  $|R^c|\geqslant2k+1$.  In the following, we further prove the following claim.
\begin{claim}\label{claim5.6b}
  $G^{*}[R^{c}]\in\mathcal{L}_{|R^{c}|,k-1}$.
\end{claim}
\begin{proof}[\bf Proof of Claim~\ref{claim5.6b}]
    If $G^*[R^c]\in \mathcal J_1\setminus  (\{S_{|R^{c}|,k-1}^{+}\}\cup  \mathcal{L}_{|R^{c}|,k-1})$, then by Lemma~\ref{lem3.3a}, we have $\eta(R)<\gamma$, a contradiction. Now we consider  $G^{*} [R^{c}]\cong S_{|R^{c}|,k-1}^{+}$. Let $R_1=\{u_1,\ldots,u_{|R_1|}\}$ be the set of dominating vertices in  $S_{|R^c|,k-1}^+$, and $R_2=\{v_1,\ldots,v_{|R_2|}\}$ be the set $R^c\setminus R_1$. It is clear that $|R_1|=k-1$ and $|R_2|\geqslant k+2$. Moreover, let $v_1v_2$ be the unique edge within $R_2$.
    Note that $d_{R^c}(u)=k$ for $u\in \{v_1,v_2\}$, $d_{R^c}(u)=k-1$ for $u\in R_2\setminus\{v_1,v_2\}$ and $d_{R^c}(u)=|R^c|-1$ for $u\in R_1$. Since $\gamma=-\frac{k(k-1)}{2} $ and $e(R^c)={k-1\choose 2}+(k-1)(|R^c|-k+1)+1$, by \eqref{4},  we obtain 
\begin{align}\label{16b}
\eta(R^c)
&=\sum_{u\in R_1}(d_{R^c}(u)-k+1)x_u+\sum_{u\in\{v_1,v_2\}}(d_{R^c}(u)-k+1)x_u-e(R^c)\notag \\
&=  (k-1)(|R^c|-k)- (|R^c|-k)\sum_{u\in R_1}(1-x_u)+x_{v_1}+x_{v_2}-e(R^c)\notag\\
    &= \gamma-1+x_{v_1}+x_{v_2}-(|R^c|-k)(k-1-\sum_{u\in R_1}x_u).
\end{align}

If $x_{v_1}+x_{v_2}<1$, then by \eqref{16b}, $\eta(R^c)<\gamma$, and so by Lemma~\ref{lem3.2a}, $\eta(R)<\gamma$, a contradiction to \eqref{6b}. In the following, we consider $x_{v_1}+x_{v_2} \geqslant 1$. 

Note that both $v_1$ and $v_2$ have no neighbor in $R \setminus R^c$. Otherwise, $N_{G^*}(u^*)$ contains a path of order $2k+1$, and so $G^*$ contains an $F_{2k+2}$, a contradiction. Therefore, $\lambda (x_{v_1}+x_{v_2})=x_{v_1}+x_{v_2}+2x_{u^*} + 2\sum_{u\in R_1}x_u + \sum_{w\in N_S(\{v_1,v_2\})}d_{\{v_1,v_2\}}(w)x_w$, and so $e(S,\{v_1,v_2\})\geqslant(\lambda-1)(x_{v_1}+x_{v_2})-2k$. 

Let $w\in N_S(\{v_1,v_2\})$. Clearly, $N_{R_1}(w) =\emptyset$. Otherwise, $G^*$ contains an $F_{2k+2}$ whose central vertex is in $N_{R_1}(w)$, a contradiction. That is to say, $d_{R_1}(w)=0$ for all $w\in N_S(\{v_1,v_2\})$. Then for all $w\in N_S(\{v_1, v_2\})$ with $d_S(w)=0,$ one has
$\lambda x_w\leqslant\lambda x_{u^*}-\sum_{u\in R_1}x_u=\lambda -\sum_{u\in R_1}x_u$, and so $x_w\leqslant 1-\frac{1}{\lambda}\sum_{u\in R_1}x_u$. 
By \eqref{6b} and~\eqref{16b}, one has $e(S)\leqslant 1$. 
Combining  with $e(S,\{v_1,v_2\})\geqslant (\lambda-1)(x_{v_1}+x_{v_2})-2k$ and $x_{v_1}+x_{v_2}\geqslant 1$ yields
$e(\hat{S},\{v_1,v_2\})\geqslant (\lambda-1)(x_{v_1}+x_{v_2})-2k-4\geqslant(\lambda-2k-5)(x_{v_1}+x_{v_2})$, where $\hat{S}=\{w\in S| d_S(w)=0\}$.

Now by~\eqref{5b} and~\eqref{6b}, 
$\eta(R)\geqslant \gamma+\frac{\lambda-2k-5}{\lambda}(x_{v_1}+x_{v_2})\sum_{u\in R_1}x_u$. Combining with \eqref{16b}, one has 
\begin{align*}
    \frac{\lambda-2k-5}{\lambda}(x_{v_1}+x_{v_2})\sum_{u\in R_1}x_u\leqslant -1+x_{v_1}+x_{v_2}-(|R^c|-k)(k-1-\sum_{u\in R_1}x_u),
\end{align*}
and so 
\begin{align}\label{3.7b}
    \left(\frac{\lambda-2k-5}{\lambda}\sum_{u\in R_1}x_u-1\right)(x_{v_1}+x_{v_2})\leqslant -(|R^c|-k)(k-1-\sum_{u\in R_1}x_u)-1.
\end{align}

If $\sum_{u\in R_1}x_u\leqslant k-\frac{3}{2}$, then by \eqref{3.7b},  we obtain $-2<(\frac{\lambda-2k-5}{\lambda}\sum_{u\in R_1}x_u-1)(x_{v_1}+x_{v_2})\leqslant-\frac{1}{2}(|R^c|-k)-1\leqslant -\frac{k+3}{2}$, a contradiction. 
If  $ k-\frac{3}{2}<\sum_{u\in R_1}x_u\leqslant k-1 $, then by \eqref{3.7b}, we obtain  $(\frac{\lambda-2k-5}{\lambda}\sum_{u\in R_1}x_u-1)(x_{v_1}+x_{v_2}) < 0$, and so $\frac{\lambda-2k-5}{\lambda}(k-\frac{3}{2})-1 < 0$,  a contradiction. Therefore, $G^*[R^c]\in \mathcal L_{|R^c|,k-1}$.

This completes the proof.
\end{proof}

Now, we come back to show Theorem~\ref{thm2.1}. 
Note that $V(G)=\{u^*\}\cup R\cup S$. It suffices to show $S=\emptyset$ and $G^*[R]\cong S_{|R|,k-1}$.

By Claim~\ref{claim5.6b}, we have $G^*[R^c]\in\mathcal L_{|R^c|,k-1}$. By Lemmas~\ref{lem3.2a} and~\ref{lem3.3a},  we get $\eta(R)\leqslant \eta(R^c)\leqslant \gamma$. By~\eqref{6b},  we know that $\eta(R)\geqslant e(S)+\gamma$. Thus, $e(S)=0$ and $\eta(R)=\eta(R^c)=\gamma$. 
 By Lemmas~\ref{lem3.2a} and~\ref{lem3.3a},  one has $R=R^c$, and $x_u=1$ for each $u\in R$ with $d_R(u)\geqslant k$. Also, combining  with  $\eta(R)=\gamma$ and \eqref{5b}, \eqref{6b}, we  have $x_w=1$ for each $w\in S$.

Next we show $S=\emptyset$. Otherwise, let $w_0$ be in $S$. Then $x_{w_0}=1=x_{u^*}$. Since $e(S)=0$, we have $N(w_0)\subseteq N(u^*)$, and so  $N(w_0)=N(u^*)$. Note that $G^*[R]\in\mathcal L_{|R|,k-1}$, that is, $G^*[R]$ is obtained from $S_{|R|,k-1}^+$ by deleting some edge $e^*$. One may still partition $R$ into $R_1\cup R_2$, where $R_1=\{u_1,\ldots,u_{|R_1|}\}$ and $R_2=
\{v_1,\ldots,v_{|R_2|}\}$.  

$\bullet$ $k=3$ and $e^*=u_1u_2$. If $|S|=1,$ then $S=\{w_0\}$ and $e^*=u_1u_2$. So we have  $\lambda =\lambda x_{v_1}\leqslant 5$, which contradicts $\lambda >\sqrt{m}\geqslant 100$. If $|S|\geqslant 2$, then there is an $F_{8}$ in $G^*$, a contradiction.

$\bullet$ $k=3,e^*\not=u_1u_2$, or $k\geqslant 4$. No matter how $e^*$ is chosen, there always exists a vertex subset $R_2'\subseteq R_2$ of size $(k+2)$ such that $G^*[\{u^*,w_0\}\cup R_1\cup R_2']$ contains a spanning subgraph being isomorphic to $F_{2k+2}$, a contradiction.
   

Now, we show $G^*[R]\cong S_{|R|,k-1}$, and more precisely, we show that $e^*=v_1v_2$.  Suppose to the contrary that $e^*\neq v_1v_2$. Then there must exist an $i\in \{1,2\}$  such that $d_R(v_i)=k$. Without loss of generality, we may assume $d_R(v_1)=k$. Recall that $x_u=1=x_{u^*}$ for each $u\in R $ with $d_R(u)\geqslant k$. Consequently, $N_{G^*}[u]=N_{G^*}[u^*]$ for each $u\in R$ with $d_R(u)\geqslant k$. Hence, $N_{G^*}[v_1]=N_{G^*}[u^*]$, which implies $v_1$ is adjacent to each vertex of $R_2$, a contradiction to the choice of $R_2$. Therefore, $e^*=v_1v_2$ and so $G^*[R]\cong S_{|R|,k-1}$.

This completes the proof.
\end{proof}


\begin{proof}[\bf Proof of Corollary~\ref{thm1.3}]
Recall that  $k\geqslant 3$ and $m\geqslant \frac{9}{4}k^6+6k^5+46k^4+56k^3+196k^2$.
It is clear that $\mathcal G(m,F_{k,3}) \subseteq  \mathcal G(m,F_{2k+2})$ and there is no  $F_{k,3}$ in $G \cong K_k\vee(\frac{m}{k}-\frac{k-1}{2})K_1$. By Theorem~\ref{thm2.1}, we obtain that  if $G\in\mathcal G(m,F_{k,3})$, 
then
$
\lambda(G)\leqslant\frac{k-1+\sqrt{4m-k^2+1}}{2}
$
with equality if and only if $G \cong K_k\vee(\frac{m}{k}-\frac{k-1}{2})K_1$. 
This completes the proof.
\end{proof}


\begin{proof}[\bf Proof of Corollary~\ref{thm2.2}]
Recall that  $k\geqslant 3$ and $m\geqslant \frac{9}{4}k^6+6k^5+46k^4+56k^3+196k^2$.
When $q\geqslant p\geqslant 3, s\geqslant r\geqslant 3,  p+q=2k+1$ and $r+s=2k+2,$ it is clear that $\mathcal G(m,\theta_{1,p,q})\cup  \mathcal G(m,\theta_{1,s,t}) \subseteq  \mathcal G(m,F_{2k+2})$ and there is no  $\theta_{1,p,q}$ or $\theta_{1,s,t}$ in $G \cong K_k\vee(\frac{m}{k}-\frac{k-1}{2})K_1$. By Theorem~\ref{thm2.1}, we obtain that  if $G\in\mathcal G(m, \theta_{1,p,q})\cup\mathcal G(m, \theta_{1,r,s})$, 
then
$
\lambda(G)\leqslant\frac{k-1+\sqrt{4m-k^2+1}}{2}
$
with equality if and only if $G \cong K_k\vee(\frac{m}{k}-\frac{k-1}{2})K_1$. 
This completes the proof.
\end{proof}

\section{\normalsize Concluding remarks}\label{sec-5}

In this paper, we focus on spectral extrema of graphs with a fixed size. First, we characterize the \( F_{2k + 2} \)-free graph of size \( m \) that attains the largest spectral radius (see Theorem \ref{thm2.1}). Second, we determine the \( F_{k, 3} \)-free graph of size \( m \) with the maximum spectral radius, which confirms Conjecture \ref{cj1-5}.  
Notably, the inclusions \( C_{2k + 1} \subseteq \theta_{1, 2, 2k - 1} \subseteq F_{2k + 1} \), \( C_{2k + 2} \subseteq \theta_{1, 2, 2k} \subseteq F_{2k + 2} \), and \( F_{2k + 1} \subseteq F_{2k + 2} \) hold. Thus, Theorem \ref{cj1-3} follows from Theorem \ref{cj1-4}, and for sufficiently large \( m \), Theorem \ref{cj1-4} can also be derived from Conjecture \ref{conj-6}. For \( k \geq 3 \), Theorem \ref{thm2.1} is stronger than Conjecture \ref{conj-6}.  Consequently, the main results (\cite[Theorems~1.3, 1.4, 1.5 and 1.6]{C01}) are direct consequences of Theorem \ref{cj1-5} when \( m \geq \frac{9}{4}k^6 + 6k^5 + 46k^4 + 56k^3 + 196k^2 \). 

It is noteworthy that we identify the \( \theta_{1,p,q} \)-free graph of size \( m \) with the largest spectral radius, where $q\geqslant p\geqslant 3$ and $p+q\geqslant 7$, thereby resolving Problem \ref{pb-2} for these cases. When combined with \cite{GL-2024,N-21,C01,LLL24,SLW-23}, the graphs in \( \mathcal{G}(m, \theta_{1,p,q}) \) that attain the maximum spectral radius have been fully determined for \( q \geqslant p \geqslant 2 \). Thus, we propose the following natural and interesting problems.  
\begin{pb}\label{pb5-3}
How can we characterize the graphs among $\mathcal{G}(m,\theta_{r,p,q})$ having the largest spectral radius for $q\geqslant p\geqslant r\geqslant 2?$
\end{pb}
Recently, Gao and Li~\cite{GL-2025} solved Problem~\ref{pb5-3} for $r=2,p=2$ and $q=3.$

\section*{\normalsize Disclosure statement}
The authors did not report any potential conflict of interest.

\subsection*{Acknowledgements}
We take this opportunity to thank the anonymous referees for their careful reading of the manuscript and suggestions which have immensely helped us in getting the article to its present form. Shuchao Li financially supported by the National Natural Science Foundation of China (Grant Nos. 12171190, 11671164), the Special Fund for Basic Scientific Research of Central Colleges (Grant Nos. CCNU25JC006, CCNU25HD044, CCNU25JCPT031) and the Open Research Fund of Key Laboratory of Nonlinear Analysis \& Applications (CCNU), Ministry of Education of China (Grant No. NAA2025ORG010). 

\section*{\normalsize Data availability}
Data sharing is not applicable to this article, as no data sets were generated or analyzed during the current study.

\end{document}